\pdfoutput=1
\documentclass[hidelinks,onefignum,onetabnum]{siamart251216}

\usepackage{amsfonts}
\usepackage{graphicx}
\usepackage{epstopdf}
\ifpdf
  \DeclareGraphicsExtensions{.eps,.pdf,.png,.jpg}
\else
  \DeclareGraphicsExtensions{.eps}
\fi

\usepackage{amsmath,amssymb,booktabs}

\usepackage{enumitem}
\setlist[itemize]{leftmargin=*}

\newcommand{\R}{\mathbb R}
\newcommand{\Sph}{\mathbb S}
\newcommand{\tr}{{\rm tr}}
\newcommand{\br}{{\rm br}}
\newcommand{\EG}{{\rm EG}}
\newcommand{\jump}[1]{[\![#1]\!]}
\newcommand{\avg}[1]{\{\!\!\{#1\}\!\!\}}

\newsiamremark{remark}{Remark}
\newsiamremark{assumption}{Assumption}
\crefname{assumption}{Assumption}{Assumptions}
\Crefname{assumption}{Assumption}{Assumptions}

\headers{Deterministic operator networks}{Young Ju Lee and Youngkyu Lee}

\title{Deterministic Operator Networks: Approximation Rates and Locally Conservative Operator Learning\thanks{Submitted to the editors \today.
\funding{Young Ju Lee was supported in part by NSF-DMS 2208499 and Youngkyu Lee was supported by a grant from Kyung Hee University in 2026 (KHU-20263365).}}}

\author{Young Ju Lee\thanks{Mathematics Department, Texas State University, San Marcos, TX, US
  (\email{yjlee@txstate.edu}).}
\and Youngkyu Lee\thanks{Department of Mathematics, Kyung Hee University, Seoul, Republic of Korea
  (\email{youngkyu\_lee@khu.ac.kr}).}}

\usepackage{amsopn}

\DeclareMathOperator{\spn}{span}

\crefname{appendix}{Appendix}{Appendices}
\Crefname{appendix}{Appendix}{Appendices}

\ifpdf
\hypersetup{
  pdftitle={Deterministic Operator Networks: Approximation Rates and Locally Conservative Operator Learning},
  pdfauthor={Young Ju Lee and Youngkyu Lee}
}
\fi

\begin{document}

\maketitle

\begin{abstract}
We introduce deterministic operator networks (DetONets), a DeepONet architecture with fixed, quasi-uniformly distributed $\mathrm{ReLU}^k$ features and a linear least-squares readout, avoiding nonconvex training of hidden parameters. An explicit construction based on orthogonal projections yields convergence rates in trunk and branch neuron counts with stability constant one, while separating the input-encoding error, which is bounded below by an irreducible floor, from the network approximation error. We further develop a constraint layer that enforces local conservation exactly for every input. For Darcy flow, we establish neural-resolution-independent stability of this layer on a fixed auxiliary mesh, including mixed Dirichlet--Neumann conditions under the stated stabilization assumptions. The framework also provides computable error certificates and supports training without solution data. Numerical experiments demonstrate the predicted convergence behavior and elementwise conservation.
\end{abstract}



\section{Introduction}\label{sec:intro}

Operator learning approximates the solution maps of parametric PDEs from input-output data, enabling predictions for unseen coefficients without a new PDE solve~\cite{bhattacharya2021model,kovachki2023neural,li2021fourier,lu2021learning}. DeepONets~\cite{lu2021learning} represent solutions using trunk functions of space and branch functions of encoded inputs. Existing error analyses separate the encoding, approximation, and reconstruction errors~\cite{cao2026operator,chen2026neural,kovachki2021universal,lanthaler2022error}; however, nonconvex training need not realize the approximant constructed in the analysis. Moreover, penalizing PDE residuals does not generally enforce exact cellwise conservation.

Extreme learning machines fix hidden parameters and train only a linear readout~\cite{dong2021local,dwivedi2020physics,huang2006extreme,lee2025nonoverlapping}. Related approaches include random features~\cite{rahimi2007random}, RandONets~\cite{fabiani2025randonets}, and deterministic linearized neural networks~\cite{kim2026locally,liu2025integral}, whose parameters are quasi-uniform on spheres. The latter admit optimal $\mathrm{ReLU}^k$ approximation rates~\cite{liu2025integral} and have been enriched to enforce local conservation for elliptic problems~\cite{kim2026locally}.

We propose a \emph{deterministic operator network} (DetONet) with fixed, quasi-uniform $\mathrm{ReLU}^k$ trunk and branch features and a linear least-squares readout. Continuous orthogonal projections provide the approximation analysis; transferring it to empirical least squares requires sampling assumptions~\cite{cohen2013stability,cohen2017optimal}. An affine correction layer then enforces discrete constraints without nonlinear training. Applied to ENN~\cite{kim2026locally}, it yields elementwise conservative predictions even for unseen inputs.
The main outcomes are as follows:
\begin{itemize}
\item \emph{Provable, deterministic operator learning.} A linear least-squares readout replaces the nonlinear training of hidden parameters, and an explicit projection construction gives trunk and branch approximation rates with stability constant one and quantifies the input-encoding error.
\item \emph{Exact conservation with controlled stability.} A constraint layer makes every predicted Darcy flux locally conservative, including for unseen inputs. On a fixed auxiliary mesh, its stability bound is independent of neural width and resolution, including mixed Dirichlet--Neumann conditions under the stated stabilization assumptions.
\item \emph{Verifiable predictions without solution data.} A computable residual provides a two-sided error certificate and enables solution-free training. The analysis also quantifies Darcy input-encoding errors for bounded conductivities and $H^{-1}$ sources; experiments confirm the predicted convergence and elementwise conservation.
\end{itemize}

Sections~\ref{sec:detonet}--\ref{sec:analysis} develop DetONet and its error bounds; \cref{sec:constrained} establishes constrained predictions and Darcy conservation, and \cref{sec:numerics} reports experiments.

Throughout this paper, $\Omega\subset\R^d$ is a bounded Lipschitz domain; for any integer $m\ge 1$, $\Sph^m$ denotes the unit sphere in $\R^{m+1}$; $\sigma_k(y)=\max(y,0)^k$ is the $\mathrm{ReLU}^k$ activation function, and
\[
r_{\max}(m,k):=\frac{m+2k+1}{2}.
\]
We write $a\lesssim b$ if $a\le Cb$ with a constant $C$ independent of the numbers of neurons.
For a set $\Xi\subset\R^q$ of positive measure, $L^2(\Xi)$ carries the normalized measure, i.e., $\|a\|_{L^2(\Xi)}^2=|\Xi|^{-1}\int_\Xi|a|^2\,dz$.

\section{Deterministic operator networks}\label{sec:detonet}

In this section, we first describe the operator learning problem together with the input encoder, and then introduce the DetONet.
Finally, we derive an error decomposition that is used throughout this paper.

\subsection{Operator learning with an input encoder}\label{sec:encoder-def}
Let $X$ be a Banach space of functions on $\Omega$, let $\mathcal U\subset X$ be a set of admissible inputs, and let $\mathcal G:\mathcal U\to H^t(\Omega)$, $t\ge0$, be the operator to be learned.
As a model problem, we consider the Darcy flow
\begin{equation}\label{eq:darcy}
\left\{\begin{aligned}
-\nabla\cdot(\kappa\nabla p)&=f&&\text{in }\Omega,\\
p&=p_D&&\text{on }\Gamma_D,\\
-\kappa\nabla p\cdot n&=j_N&&\text{on }\Gamma_N,
\end{aligned}\right.
\end{equation}
where $\partial\Omega=\Gamma_D \cup \Gamma_N$ with $|\Gamma_D|>0$, and the boundary data $p_D$ and $j_N$ are fixed.
The input is either the permeability $u=\kappa$ or the pair $u=(\kappa,f)$ of the permeability and the source, and $\mathcal G(u)=p$ is the pressure.

We assume a finite-dimensional input encoder
\[
E:\mathcal U\to\R^q,\qquad u\mapsto z=E(u).
\]
Examples include sensor values followed by dimensionality reduction, projection coefficients, and principal components.
Let $\Xi\subset\R^q$ be a bounded Lipschitz domain with $E(\mathcal U)\subset\Xi$.

The network approximates a target $g:\Xi\to H^t(\Omega)$ through $\mathcal G(u)\approx g(E(u))$; write $g(z,x)=g(z)(x)$. The target is used only in the analysis, not evaluated during training.
We construct $g$ from a \emph{reconstruction map}
\begin{equation}\label{eq:recon}
R:\Xi\to\mathcal U\qquad\text{with}\qquad E(R(z))=z\quad\text{for }z\in E(\mathcal U),
\end{equation}
which rebuilds an admissible input from its code, and set
\begin{equation}\label{eq:g}
g:=\mathcal G\circ R,\qquad u_q:=R(E(u)).
\end{equation}
Then $g(E(u))=\mathcal G(u_q)$ is the exact solution for the rebuilt input $u_q$, and the condition $E\circ R=\mathrm{id}$ means that $u$ and $u_q$ have the same code.
Here, we assume $R(\Xi)\subset\mathcal U$, i.e., $\mathcal U$ contains the reconstructions of all $z\in\Xi$, not only of the codes of admissible inputs.
The following two examples are central in this paper. 
\begin{itemize}
\item \emph{Exact encoding.} If $E$ is a bijection from $\mathcal U$ onto $\Xi$, then $R=E^{-1}$, $g=\mathcal G\circ E^{-1}$, and $\mathcal G=g\circ E$ holds exactly.
\item \emph{Darcy flow.} Let $E(\kappa)\in\R^q$ be the coordinates of $\log\kappa$ in a basis $e_1,\dots,e_q$, and set
\begin{equation}\label{eq:kappa-z}
R(z):=\kappa_z:=\exp\Big(\sum_{\ell=1}^qz_\ell e_\ell\Big),\qquad g(z)=p(\kappa_z).
\end{equation}
Then $\kappa_q:=\kappa_{E(\kappa)}$ is the permeability $\kappa$ with $\log\kappa$ reduced to $q$ terms.
Since $\kappa_z$ is bounded above and below for $z$ in the bounded set $\Xi$, $p(\kappa_z)$ is well defined, and $g$ is holomorphic in $z$, as is well known for parametric elliptic problems~\cite{cohen2015approximation}.
With a source, i.e., $u=(\kappa,f)$, we use $E(u)=(z,\zeta)\in\R^{q_\kappa}\times\R^{q_f}$ with $q=q_\kappa+q_f$ and $R(z,\zeta)=(\kappa_z,f_\zeta)$ with $f_\zeta=\sum_{j=1}^{q_f}\zeta_j\varphi_j$ for a basis $\varphi_1,\dots,\varphi_{q_f}$; then $g$ is affine in $\zeta$, since the boundary data are fixed.
\end{itemize}

\subsection{The deterministic operator network}\label{sec:detonet-def}
The trunk and branch features of the DetONet are given by
\begin{align*}
T_j(x)&:=\sigma_{k_\tr}\big(\theta^\tr_j\cdot(x,1)\big),\quad j=1,\dots,N,\\
B_i(z)&:=\sigma_{k_\br}\big(\theta^\br_i\cdot(z,1)\big),\quad i=1,\dots,M,
\end{align*}
where the trunk parameters $\theta^\tr_j$ are fixed and quasi-uniformly distributed on $\Sph^d$, and the branch parameters $\theta^\br_i$ are fixed and quasi-uniformly distributed on $\Sph^q$.
With $T(x)=(T_1(x),\dots,T_N(x))^\top$, $B(z)=(B_1(z),\dots,B_M(z))^\top$, and a readout $W=(w_{ji})\in\R^{N\times M}$, the DetONet is defined as
\begin{equation}\label{eq:FMN}
F_{M,N}[W](z,x):=\sum_{j=1}^N\sum_{i=1}^M w_{ji}\,T_j(x)\,B_i(z)=T(x)^\top W\,B(z),
\end{equation}
and its prediction for the input $u$ is $F_{M,N}[W](E(u),\cdot)$.
Since only the readout is trained, the DetONets form the linear space
\[
\mathcal V_{M,N}:=\spn\{B_i\otimes T_j\}=\{F_{M,N}[W]:\ W\in\R^{N\times M}\},\qquad\dim\mathcal V_{M,N}\le NM.
\]
Given training inputs $u_1,\dots,u_{n_s}$ with codes $z_m=E(u_m)$ and outputs evaluated at grid points $x_1,\dots,x_{n_x}$, let $\mathsf T_{lj}=T_j(x_l)$, $\mathsf B_{im}=B_i(z_m)$, and $\mathsf V_{lm}=\mathcal G(u_m)(x_l)$.
Then the readout is computed by the least squares solution
\begin{equation}\label{eq:ls}
W=\mathsf T^\dagger\,\mathsf V\,\mathsf B^\dagger,
\end{equation}
where ${}^\dagger$ denotes the Moore--Penrose inverse.
Note that no gradient-based optimization is required, and that the DetONet is a DeepONet with shallow trunk and branch networks whose inner parameters are fixed in advance.

\subsection{Error decomposition}\label{sec:errid}
For $u\in\mathcal U$ and $z=E(u)$, we define the \emph{modeling error} of $E$ and the \emph{network error} as
\[
\mathrm{err}_E(u):=\mathcal G(u)-g(E(u)),\qquad \mathrm{err}_N(z):=g(z)-F_{M,N}[W](z,\cdot).
\]
Adding and subtracting $g(E(u))$ gives the identity
\begin{equation}\label{eq:identity}
\mathcal G(u)-F_{M,N}[W](E(u),\cdot)=\mathrm{err}_E(u)+\mathrm{err}_N(E(u)).
\end{equation}
The two terms are of different nature.
The modeling error depends on the encoder and target, not on neural width or training, and it vanishes for exact encoding.
We measure it by
\[
\varepsilon_E:=\sup_{u\in\mathcal U}\|\mathcal G(u)-g(E(u))\|_{H^t(\Omega)}.
\]
Since inputs with the same code receive the same prediction, no choice of $g$, and hence no network, can do better than the spread of their solutions.
That is, for every $g$,
\begin{equation}\label{eq:floor}
\varepsilon_E\ \ge\ \varepsilon_E^\ast:=\sup_{z\in E(\mathcal U)}\ \inf_{v\in H^t(\Omega)}\ \sup_{u\in\mathcal U:\,E(u)=z}\|\mathcal G(u)-v\|_{H^t(\Omega)}.
\end{equation}
The floor $\varepsilon_E^\ast$ depends only on $E$; a Chebyshev center on each encoder fiber may attain it but need not give a smooth target, whereas the reconstruction~\eqref{eq:kappa-z} yields a holomorphic one.
The network error depends on the regularity of $g$ in $(z,x)$; we estimate the two terms in \cref{sec:network,sec:model}.

\section{Error analysis}\label{sec:analysis}

In this section, we estimate the two terms of the identity~\eqref{eq:identity}.
We first bound the network error by constructing an explicit DetONet, and then bound the modeling error by the stability of the solution operator of the Darcy problem.

\subsection{Network error}\label{sec:network}
Approximation is measured by a linear quasi-interpolant, and stability is provided by orthogonal projections in each variable.

Let us recall the simultaneous approximant of~\cite[Lemma~A.1 and the proof of Proposition~2.3]{kim2026locally}, which is built on the integral representations of Liu, Mao, and Xu~\cite{liu2025integral}.
Let $\{\theta_j\}_{j=1}^n\subset\Sph^d$ be quasi-uniform with mesh size $\delta=n^{-1/d}$.
Then there is a linear operator $I:H^r(\Omega)\to\spn\{\sigma_k(\theta_j\cdot(\,\cdot\,,1))\}_{j=1}^n$ such that
\begin{equation}\label{eq:qi}
\|v-Iv\|_{H^t(\Omega)}\lesssim\delta^{\,r-t}\,\|v\|_{H^r(\Omega)},\qquad 0\le t\le\min\{k,r\},\quad 0<r\le r_{\max}(d,k).
\end{equation}
At saturation, i.e., $r=r_{\max}(d,k)$, the rate is $n^{-\frac12-\frac{2(k-t)+1}{2d}}$.
The same holds on $\Xi\subset\R^q$ with $d$ replaced by $q$.
Note that we use $I$ only to bound best approximation errors, so that $I$ need not be bounded on $H^t(\Omega)$ uniformly in $n$.

We first define the trunk projection.
Fix $0\le t\le k_\tr$, so that $T_j\in H^t(\Omega)$, and let $\mathcal S_N:=\spn\{T_1,\dots,T_N\}\subset H^t(\Omega)$.
Let $P_\tr:H^t(\Omega)\to\mathcal S_N$ be the orthogonal projection in the $H^t(\Omega)$ inner product, and let $I_\tr$ be the quasi-interpolant~\eqref{eq:qi} with $k=k_\tr$ and the trunk parameters.

\begin{lemma}[Trunk projection]\label{lem:trunk}
For $v\in H^t(\Omega)$, it holds that $\|P_\tr v\|_{H^t(\Omega)}\le\|v\|_{H^t(\Omega)}$.
For $v\in H^r(\Omega)$ with $t\le\min\{k_\tr,r\}$ and $0<r\le r_{\max}(d,k_\tr)$, we have
\[
\|v-P_\tr v\|_{H^t(\Omega)}\le\|v-I_\tr v\|_{H^t(\Omega)}\lesssim N^{-(r-t)/d}\,\|v\|_{H^r(\Omega)}.
\]
\end{lemma}

\begin{proof}
An orthogonal projection has norm at most one.
Since $I_\tr v\in\mathcal S_N$ and $P_\tr v$ is the best approximation from $\mathcal S_N$ in $H^t(\Omega)$, we obtain $\|v-P_\tr v\|_{H^t}\le\|v-I_\tr v\|_{H^t}$.
Then the result follows from~\eqref{eq:qi} with $\delta=N^{-1/d}$.
This completes the proof.
\end{proof}

Note that the stability constant is one for every $N$, independently of the conditioning of the features $T_j$, because $P_\tr$ is defined through the space $\mathcal S_N$ rather than through the coefficients in the spanning set.

We next define the branch projection.
Let $\mathcal B_M:=\spn\{B_1,\dots,B_M\}\subset L^2(\Xi)$, let $P_\br:L^2(\Xi)\to\mathcal B_M$ be the orthogonal projection, and let $I_\br$ be the quasi-interpolant~\eqref{eq:qi} on $\Xi$ with $k=k_\br$ and the branch parameters.
For a separable Hilbert space $H$, the space $\mathcal B_M\otimes H=\{\sum_ih_iB_i:\ h_i\in H\}$ is a closed subspace of $L^2(\Xi;H)$.
We let $P_\br$ and $I_\br$ act on $a:\Xi\to H$ componentwise.
That is, for an orthonormal basis $\{\phi_m\}$ of $H$ and $a_m:=(a,\phi_m)_H$, we set $(P_\br a,\phi_m)_H=P_\br a_m$, and likewise for $I_\br$.

\begin{lemma}[Hilbert-valued branch projection]\label{lem:branch}
For $a\in L^2(\Xi;H)$, it holds that $\|P_\br a\|_{L^2(\Xi;H)}\le\|a\|_{L^2(\Xi;H)}$.
For $0<s\le r_{\max}(q,k_\br)$ and $a\in H^s(\Xi;H)$, we have
\[
\|a-P_\br a\|_{L^2(\Xi;H)}\le\|a-I_\br a\|_{L^2(\Xi;H)}\lesssim M^{-s/q}\,\|a\|_{H^s(\Xi;H)}
\]
with the constant of the scalar estimate~\eqref{eq:qi}.
\end{lemma}

\begin{proof}
By Parseval's identity, $\|a-P_\br a\|^2_{L^2(\Xi;H)}=\sum_m\|a_m-P_\br a_m\|^2_{L^2(\Xi)}$, so that $P_\br$ is the orthogonal projection of $L^2(\Xi;H)$ onto $\mathcal B_M\otimes H$.
Hence, it has norm at most one and gives the best approximation from $\mathcal B_M\otimes H$, which contains $I_\br a$; this proves the first inequality.
For the second inequality, Parseval's identity and~\eqref{eq:qi} on $\Xi$ with $t=0$, $r=s$, and $\delta=M^{-1/q}$ give
\[
\|a-I_\br a\|_{L^2(\Xi;H)}^2=\sum_m\|a_m-I_\br a_m\|_{L^2(\Xi)}^2\lesssim\delta^{2s}\sum_m\|a_m\|_{H^s(\Xi)}^2=\delta^{2s}\|a\|_{H^s(\Xi;H)}^2.
\]
This completes the proof.
\end{proof}

In order to combine the trunk and branch projections, let $P^x:=\mathrm{id}_z\otimes P_\tr$, i.e., $(P^xg)(z,\cdot)=P_\tr\,g(z,\cdot)$, and $P^z:=P_\br\otimes\mathrm{id}_x$.
Both are orthogonal projections of $L^2(\Xi;H^t(\Omega))$ acting on different variables.
Hence, they commute, and $P^xP^z$ is the orthogonal projection of $L^2(\Xi;H^t(\Omega))$ onto $\mathcal B_M\otimes\mathcal S_N=\mathcal V_{M,N}$.
We now take for the DetONet $F_{M,N}[W^\ast]:=P^xP^zg$, where $W^\ast$ are the coefficients of $P^xP^zg\in\mathcal V_{M,N}$ in the spanning set $\{B_iT_j\}$; if the features are linearly dependent, any choice gives the same function.
Then $F_{M,N}[W^\ast]$ is the best approximation of $g$ from $\mathcal V_{M,N}$ in $L^2(\Xi;H^t(\Omega))$.
We split its network error as
\[
g-F_{M,N}[W^\ast]=\mathrm{err}_\tr+\mathrm{err}_\br,\qquad\mathrm{err}_\tr:=g-P^xg,\qquad\mathrm{err}_\br:=P^x(g-P^zg).
\]
Note that the two terms are orthogonal in $L^2(\Xi;H^t(\Omega))$, since $\mathrm{err}_\tr$ lies in the range of $\mathrm{id}-P^x$ and $\mathrm{err}_\br$ in the range of $P^x$.

\begin{theorem}[Approximation rate of the DetONet]\label{thm:rate}
Let $0\le t\le k_\tr$, and suppose that
\[
g\in L^2\big(\Xi;H^{r_{\max}(d,k_\tr)}(\Omega)\big)\cap H^{r_{\max}(q,k_\br)}\big(\Xi;H^t(\Omega)\big).
\]
Define the worst-case error at the code $z$ by
\[
e(z):=\sup\big\{\|\mathcal G(u)-F_{M,N}[W^\ast](z,\cdot)\|_{H^t(\Omega)}:\ u\in\mathcal U,\ E(u)=z\big\},
\]
with $e(z):=0$ if no such $u$ exists.
Then we have
\begin{equation}\label{eq:rate}
\begin{aligned}
\|e\|_{L^2(\Xi)}\le\varepsilon_E
&+C\,N^{-\frac12-\frac{2(k_\tr-t)+1}{2d}}\,\|g\|_{L^2(\Xi;H^{r_{\max}(d,k_\tr)}(\Omega))}\\
&+C\,M^{-\frac12-\frac{2k_\br+1}{2q}}\,\|g\|_{H^{r_{\max}(q,k_\br)}(\Xi;H^t(\Omega))},
\end{aligned}
\end{equation}
where $C$ depends only on $d$, $q$, $k_\tr$, $k_\br$, $t$, $\Omega$, $\Xi$, and the quasi-uniformity constants.
More generally, if $g\in L^2(\Xi;H^r(\Omega))\cap H^s(\Xi;H^t(\Omega))$ with $t\le\min\{k_\tr,r\}$, $r\le r_{\max}(d,k_\tr)$, and $0<s\le r_{\max}(q,k_\br)$, then the trunk and branch terms are of order $N^{-(r-t)/d}$ and $M^{-s/q}$, respectively.
\end{theorem}

Here, $e$ is assumed to be measurable; otherwise, $\|e\|_{L^2(\Xi)}$ is an upper integral.
If the codes $E(u)$ of random inputs have a density bounded by $\varrho$ with respect to the normalized measure on $\Xi$, then \eqref{eq:rate} bounds the root-mean-square error over the inputs, with the factor $\varrho^{1/2}$.

\begin{proof}
For $u\in\mathcal U$ with $E(u)=z$, the identity~\eqref{eq:identity} and the triangle inequality give $\|\mathcal G(u)-F_{M,N}[W^\ast](z,\cdot)\|_{H^t}\le\varepsilon_E+\|(\mathrm{err}_\tr+\mathrm{err}_\br)(z)\|_{H^t}$, and hence $e(z)\le\varepsilon_E+\|(\mathrm{err}_\tr+\mathrm{err}_\br)(z)\|_{H^t}$.
Note that, in the normalized measure, the constant $\varepsilon_E$ has $L^2(\Xi)$ norm $\varepsilon_E$.
Applying \cref{lem:trunk} for each $z$ and integrating over $\Xi$, we obtain
\[
\|\mathrm{err}_\tr\|_{L^2(\Xi;H^t(\Omega))}\lesssim N^{-(r-t)/d}\,\|g\|_{L^2(\Xi;H^r(\Omega))}.
\]
Since $\|P^x\|\le1$, \cref{lem:branch} with $H=H^t(\Omega)$ gives
\[
\|\mathrm{err}_\br\|_{L^2(\Xi;H^t(\Omega))}\le\|g-P^zg\|_{L^2(\Xi;H^t(\Omega))}\lesssim M^{-s/q}\,\|g\|_{H^s(\Xi;H^t(\Omega))}.
\]
Finally, we take $r=r_{\max}(d,k_\tr)$ and $s=r_{\max}(q,k_\br)$, and use $(r_{\max}(m,k)-t)/m=\frac12+\frac{2(k-t)+1}{2m}$.
This completes the proof.
\end{proof}

The branch exponent $\frac12+\frac{2k_\br+1}{2q}$ increases with $k_\br$ and decreases with the code dimension $q$, but it never falls below $\frac12$.
A bound for every $z$, rather than in $L^2(\Xi)$, requires estimating the branch term in $H^{s'}(\Xi)$ with $s'>q/2$ and using the Sobolev embedding; this lowers the exponent by $s'/q$ and requires $k_\br\ge s'$.
The rates saturate at $r_{\max}$, which is consistent with the lower bound of Mao and Xu on the sphere~\cite{mao2025sharp}, and one derivative in the output norm costs $1/d$ in the trunk exponent.
For example, for $d=q=2$ and $k_\tr=k_\br=3$, \eqref{eq:rate} predicts the trunk exponents $2.25$ in $L^2$ and $1.75$ in $H^1$, and the branch exponent $2.25$ for a smooth dependence on $z$.

\subsection{Modeling error}\label{sec:model}
Unlike the network error, the modeling error involves the PDE, and we bound it by the stability of the solution operator.
With the choice~\eqref{eq:g} of the target map, the modeling error contains no network; it is the difference of two exact solutions, one for the exact input and one for the rebuilt input:
\begin{equation}\label{eq:errE}
\mathrm{err}_E(u)=\mathcal G(u)-\mathcal G(u_q),\qquad u_q=R(E(u)).
\end{equation}
Suppose that $\mathcal G$ is Lipschitz continuous from a normed space $Y\supset\mathcal U$ to $H^t(\Omega)$ on $\mathcal U$ with constant $L_{\mathcal G}$.
Since $R(\Xi)\subset\mathcal U$, both $u$ and $u_q$ are admissible, so that~\eqref{eq:errE} and~\eqref{eq:floor} give
\begin{equation}\label{eq:model}
\|\mathrm{err}_E(u)\|_{H^t(\Omega)}\le L_{\mathcal G}\,\|u-u_q\|_Y,\qquad
\varepsilon_E^\ast\le\varepsilon_E\le L_{\mathcal G}\sup_{u\in\mathcal U}\|u-R(E(u))\|_Y.
\end{equation}
The two factors are of different nature: $L_{\mathcal G}$ is the stability constant of the PDE, which dictates the norm $Y$, while $\|u-u_q\|_Y$ measures how well $q$ numbers describe the input.
For elliptic problems, $Y$ is $L^\infty$ for the conductivity and a negative norm for the source, as we show next.

We consider the Darcy problem~\eqref{eq:darcy}.
Let $H^1_D(\Omega)=\{v\in H^1(\Omega):\ v=0\text{ on }\Gamma_D\}$ with the norm $\|\nabla v\|_{L^2(\Omega)}$, and let $\|\cdot\|_{-1}$ be the norm of its dual space $H^{-1}_D(\Omega)$.
We assume that
\begin{equation}\label{eq:kappa-bound}
\|\log\kappa\|_{L^\infty(\Omega)}\le\beta\qquad\text{for all }\kappa\in\mathcal U\cup R(\Xi),
\end{equation}
so that $e^{-\beta}\le\kappa\le e^\beta$, and let $C_{\nabla p}:=\sup\|\nabla p(\kappa,f)\|_{L^2(\Omega)}$ over all $\kappa$ satisfying~\eqref{eq:kappa-bound} and all admissible sources $f$.
For homogeneous boundary data, we have $C_{\nabla p}\le e^\beta\sup\|f\|_{-1}$.

\begin{lemma}[Stability of the Darcy problem]\label{lem:darcy}
For admissible $(\kappa,f)$ and $(\kappa',f')$, let $\varepsilon:=\|\log\kappa-\log\kappa'\|_{L^\infty(\Omega)}$.
Then it holds that
\begin{equation}\label{eq:darcy-stab}
\|\nabla(p(\kappa,f)-p(\kappa',f'))\|_{L^2(\Omega)}\le e^{2\beta}C_{\nabla p}\,\varepsilon+e^\beta\,\|f-f'\|_{-1}.
\end{equation}
Moreover, in the energy norm, the conductivity part is free of the contrast:
\begin{equation}\label{eq:darcy-sinh}
\|\kappa^{1/2}\nabla(p(\kappa,f)-p(\kappa',f))\|_{L^2(\Omega)}\le2\sinh\big(\tfrac12\varepsilon\big)\,\|\kappa'^{1/2}\nabla p(\kappa',f)\|_{L^2(\Omega)}.
\end{equation}
\end{lemma}

The proof is standard and is given in \cref{sm:sec-darcy}.
Applying \cref{lem:darcy} with $(\kappa',f')=(\kappa_q,f_q)$, we obtain
\begin{equation}\label{eq:darcy-model}
\|\nabla\mathrm{err}_E(\kappa,f)\|_{L^2(\Omega)}\le e^{2\beta}C_{\nabla p}\,\|\log\kappa-\log\kappa_q\|_{L^\infty(\Omega)}+e^\beta\,\|f-f_q\|_{-1}.
\end{equation}
That is, $Y=L^\infty$ for $\log\kappa$ and $Y=H^{-1}_D$ for $f$.
Since $2\sinh(\varepsilon/2)=\varepsilon+O(\varepsilon^3)$, the energy form~\eqref{eq:darcy-sinh} has the same first-order behavior as~\eqref{eq:darcy-stab}, but without the factor $e^{2\beta}$.

It remains to bound the encoding errors of $\log\kappa$ and $f$.
Suppose that the encoder uses point values of $\log\kappa$, i.e., $E(\kappa)=\mathbf E\,(\log\kappa)(\hat{\mathbf x})$, let $\mathbf S\in\R^{m\times q}$ with $\mathbf S_{i\ell}=e_\ell(\hat x_i)$ be the sensor matrix, and assume $\mathbf E\mathbf S=I$, e.g., $\mathbf E=\mathbf S^\dagger$ with unisolvent sensors.
Then $E\circ R=\mathrm{id}$ holds for the reconstruction~\eqref{eq:kappa-z}.
Let
\[
\Lambda:=\sup_{x\in\Omega}\sum_{i=1}^m\Big|\sum_{\ell=1}^qe_\ell(x)\,\mathbf E_{\ell i}\Big|
\]
be the Lebesgue constant of the fit, which is computable from $\mathbf E$ and $\{e_\ell\}$.

\begin{lemma}[Encoding of $\log\kappa$]\label{lem:enc}
If $\mathbf E\mathbf S=I$, then it holds that
\[
\|\log\kappa-\log\kappa_q\|_{L^\infty(\Omega)}\le(1+\Lambda)\inf_{v\in\spn\{e_\ell\}}\|\log\kappa-v\|_{L^\infty(\Omega)}.
\]
\end{lemma}

\begin{proof}
Let $Qa:=\sum_\ell(\mathbf E\,a(\hat{\mathbf x}))_\ell\,e_\ell$, so that $\log\kappa_q=Q\log\kappa$.
Since $\mathbf E\mathbf S=I$, we have $Qv=v$ for $v\in\spn\{e_\ell\}$ and $\|Qa\|_{L^\infty}\le\Lambda\|a\|_{L^\infty}$.
Hence, the result follows from $\log\kappa-Q\log\kappa=(\log\kappa-v)-Q(\log\kappa-v)$ for any $v$ in the span.
This completes the proof.
\end{proof}

In~\eqref{eq:darcy-model}, the source is truncated in $H^{-1}_D$, not in $L^2$, and this is where the negative norm matters.
Let $f_q$ be the $L^2$ projection of $f$ onto the first $q_f$ eigenfunctions of $-\Delta$ with the boundary conditions of $H^1_D(\Omega)$, whose eigenvalues satisfy $\lambda_j\sim j^{2/d}$.
Expanding $f-f_q$ in the remaining eigenfunctions gives
\begin{equation}\label{eq:negnorm}
\|f-f_q\|_{-1}\le\lambda_{q_f+1}^{-1/2}\,\|f-f_q\|_{L^2(\Omega)}\sim q_f^{-1/d}\,\|f-f_q\|_{L^2(\Omega)},
\end{equation}
so that the encoding error of the source is one order smaller than its $L^2$ error.
Note that an analysis that measures the encoding error in the $L^2$ norm of the input and transfers it to the output with a Lipschitz constant, as in~\cite{lanthaler2022error}, does not see this factor.
For the conductivity, $\|\log\kappa-\log\kappa_q\|_{L^\infty}$ decays like $1+\Lambda$ times the $L^\infty$ best approximation from $\spn\{e_\ell\}$, e.g., like $q_\kappa^{-s/d}$, up to logarithmic factors, for $\log\kappa\in C^s(\overline\Omega)$ and a boundary-adapted wavelet basis, or for a periodic $\log\kappa$ and a trigonometric basis.
Combining these bounds with \cref{thm:rate} makes every term of~\eqref{eq:rate} explicit.

\begin{corollary}[Error bound for the Darcy problem]\label{cor:darcy}
Let $t\in\{0,1\}$ and $u=(\kappa,f)$, and let $E$ and $R$ be as above with $\mathbf E\mathbf S=I$ and~\eqref{eq:kappa-bound}.
If $g=\mathcal G\circ R$ satisfies the assumptions of \cref{thm:rate}, then it holds that
\begin{multline*}
\|e\|_{L^2(\Xi)}\le C_\Omega\Big(e^{2\beta}C_{\nabla p}(1+\Lambda)\sup_{\kappa}\inf_{v\in\spn\{e_\ell\}}\|\log\kappa-v\|_{L^\infty}+e^\beta\sup_f\|f-f_q\|_{-1}\Big)\\
+\mathrm{T}_N+\mathrm{B}_M,
\end{multline*}
where $\mathrm T_N$ and $\mathrm B_M$ are the trunk and branch terms of~\eqref{eq:rate}, and $C_\Omega$ depends only on the Poincar\'e constant $C_P$ of $H^1_D(\Omega)$, i.e., $C_\Omega=C_P$ for $t=0$ and $C_\Omega=(1+C_P^2)^{1/2}$ for $t=1$.
\end{corollary}

\begin{proof}
Since $\mathrm{err}_E\in H^1_D(\Omega)$, we have $\|\mathrm{err}_E\|_{H^t}\le C_\Omega\|\nabla\mathrm{err}_E\|_{L^2}$, and the result follows from~\eqref{eq:darcy-model}, \cref{lem:enc}, and~\eqref{eq:rate}.
\end{proof}

\section{Constrained DetONets and local conservation}\label{sec:constrained}

So far, the trunk consists of $\mathrm{ReLU}^k$ neurons that approximate a function of $x$.
In many applications, however, the output should lie in a fixed discrete space, such as a finite element or ENN space, and every prediction should satisfy discrete constraints exactly, such as local conservation or boundary conditions.
In this section, we extend the DetONet to this setting, in which the trunk becomes a basis of the discrete space and a constraint layer is applied after the readout, and apply it to the ENN method for the Darcy flow~\cite{kim2026locally}.

\subsection{The constraint layer}\label{sec:layer}
For each fixed neural dimension and auxiliary partition, let $V_h$ be a finite-dimensional space, the same for all inputs $u$, let $a_u:V_h\times V_h\to\R$ be a bilinear form, $F_u\in V_h'$, and $\|\cdot\|_u$ a norm on $V_h$ such that
\begin{equation}\label{eq:coerc}
a_u(v,v)\ge\|v\|_u^2,\qquad|a_u(v,w)|\le C_a\|v\|_u\|w\|_u\qquad\forall\,v,w\in V_h,
\end{equation}
with $C_a$ independent of $u$.
The \emph{discrete solution operator} $\mathcal G_h:u\mapsto P_h(u)\in V_h$ is defined by
\[
a_u(P_h(u),w)=F_u(w)\qquad\forall\,w\in V_h.
\]
It approximates $\mathcal G$ with a discretization error $\|\mathcal G(u)-P_h(u)\|$ given by the analysis of the discrete method.
We split
\[
V_h=V_n+V_0,
\]
where $V_n=\spn\{\phi_1,\dots,\phi_{n_V}\}$ with $n_V:=\dim V_n$ is the part predicted by the network and $V_0$ is the space of the constraints.
Note that the sum need not be direct, and that the basis $\{\phi_j\}$ of $V_n$ takes the role of the trunk.

For an input $u$ and $P_n\in V_n$, we define the \emph{constraint layer} $\mathcal L_u(P_n):=P_n+P_0$, where $P_0\in V_0$ solves
\begin{equation}\label{eq:layer-gen}
a_u(P_n+P_0,\psi)=F_u(\psi)\qquad\forall\,\psi\in V_0.
\end{equation}
By~\eqref{eq:coerc}, $a_u$ is coercive on $V_0$, so that $P_0$ exists uniquely, and the layer costs one solve with the matrix of $a_u$ on $V_0$.
Note that every output of the layer satisfies the constraints~\eqref{eq:layer-gen} exactly, whatever $P_n$ is.

\begin{lemma}[Stability of the constraint layer]\label{lem:layer}
Let $g_u(v,w):=\frac12(a_u(v,w)+a_u(w,v))$ and suppose
$g_u(v,v)=\|v\|_u^2$ and $a_u=g_u$ on $V_0\times V_0$. Define the homogeneous layer
$S_uv:=v-T_av$, where $T_av\in V_0$ is uniquely determined by
$a_u(T_av,\psi)=a_u(v,\psi)$ for every $\psi\in V_0$, and set
\[
C_{\rm layer}(u):=\sup_{0\ne v\in V_n}
\frac{\|S_uv\|_u}{\|v\|_u}.
\]
Then, for any $P_n\in V_n$ and any splitting $P_h(u)=P_{h,n}+P_{h,0}$ with $P_{h,n}\in V_n$ and $P_{h,0}\in V_0$,
\[
\|P_h(u)-\mathcal L_u(P_n)\|_u
\le C_{\rm layer}(u)\|P_{h,n}-P_n\|_u.
\]
Moreover, writing $T_g:V_n\to V_0$ for the $g_u$-orthogonal projection and $D_uv\in V_0'$ for $(D_uv)(\psi):=a_u(v,\psi)-g_u(v,\psi)$, one has the exact identity
\[
\|S_uv\|_u^2=\|v-T_gv\|_u^2+\|D_uv\|_{(g_u|_{V_0})^{-1}}^2.
\]
In particular, $C_{\rm layer}(u)^2\le1+\sup_{0\ne v\in V_n}\|D_uv\|_{(g_u|_{V_0})^{-1}}^2/\|v\|_u^2$.
\end{lemma}

\begin{proof}
Since $P_h(u)$ satisfies~\eqref{eq:layer-gen}, the uniqueness of $P_0$ gives $\mathcal L_u(P_{h,n})=P_h(u)$.
Hence, $P_h(u)-\mathcal L_u(P_n)$ is the difference of two constrained outputs, which satisfies the homogeneous
constraint equations and therefore equals $S_u(P_{h,n}-P_n)$.
On $V_0$, $a_u=g_u$: this follows whenever the skew part of $a_u$
vanishes on $V_0\times V_0$, as it does for the ENN application below.
Under this condition $T_av=T_gv+(g_u|_{V_0})^{-1}D_uv$.
The two terms $v-T_gv$ and $(g_u|_{V_0})^{-1}D_uv$
are $g_u$-orthogonal. Pythagoras proves the identity, and the
orthogonal projection estimate proves the final inequality.
This completes the proof.
\end{proof}

The continuity constant of $a_u$ already bounds $C_{\rm layer}$.
Indeed, let $J_u:V_h\to V_h$ be defined by $g_u(J_uv,w)=a_u(v,w)-g_u(v,w)$ for all $w\in V_h$.
Then $J_u$ is skew-adjoint with respect to $g_u$, so that $\|v+J_uv\|_u^2=\|v\|_u^2+\|J_uv\|_u^2$, and $\|v+J_uv\|_u=\sup_{w}a_u(v,w)/\|w\|_u\le C_a\|v\|_u$ by~\eqref{eq:coerc}.
Since $(D_uv)(\psi)=g_u(J_uv,\psi)$ and $g_u(\psi,\psi)=\|\psi\|_u^2$, we obtain $\|D_uv\|_{(g_u|_{V_0})^{-1}}\le\|J_uv\|_u\le(C_a^2-1)^{1/2}\|v\|_u$, and the last bound of \cref{lem:layer} yields
\begin{equation}\label{eq:layer-Ca}
C_{\rm layer}(u)\le C_a.
\end{equation}
That is, the layer enforces the constraints at the cost of at most the factor $C_a$; for the ENN method, however, $C_a$ is not known to be bounded uniformly in the neural dimension, and we bound $C_{\rm layer}$ directly in \cref{sec:cons}.
With the branch $B(z)$ of \cref{sec:detonet-def} acting on $z=E(u)$, we define the \emph{constrained DetONet} as
\begin{equation}\label{eq:cdet}
\widehat P[W](u):=\mathcal L_u\Big(\sum_{j=1}^{n_V}\big(WB(E(u))\big)_j\,\phi_j\Big),\qquad W\in\R^{n_V\times M}.
\end{equation}
For every $W$ and every input $u$, including inputs not seen in training, $\widehat P[W](u)$ satisfies~\eqref{eq:layer-gen} exactly.
Moreover, since $\mathcal L_u$ is affine, $\widehat P[W](u)$ is affine in $W$, so that every quadratic loss in $\widehat P$ is a convex quadratic function of $W$ and training remains a linear least squares problem.

In order to estimate the error of the constrained DetONet, let $\Pi:V_h\to V_n$ be a fixed linear map with $v-\Pi v\in V_0$ for all $v\in V_h$, which exists because $V_h=V_n+V_0$.
Let $\|\cdot\|_{\rm ref}$ be a reference norm on $V_h$ with
\[
\|v\|_u\le c_+\|v\|_{\rm ref}\qquad\forall\,v\in V_h,\ u\in\mathcal U\cup R(\Xi),
\]
and let $\{\phi_j\}$ be orthonormal in $\|\cdot\|_{\rm ref}$; in practice, this is a reduced basis of $V_n$, which also removes the redundancy among the neurons.
With the reconstruction $R$ of~\eqref{eq:recon}, the target of the branch is the coefficient vector $c(z)\in\R^{n_V}$ defined by
\[
\Pi P_h(R(z))=\sum_{j=1}^{n_V}c_j(z)\,\phi_j,\qquad z\in\Xi.
\]
Let $W^\ast$ be the coefficients of its branch projection, i.e., $W^\ast B=P_\br c$ is the orthogonal projection of $c$ onto $\mathcal B_M\otimes\R^{n_V}$ in $L^2(\Xi;\R^{n_V})$.

\begin{theorem}[Error of the constrained DetONet]\label{thm:cdet}
Suppose $C_{\rm layer}(u)\le C_L$ uniformly over admissible inputs. For every input $u$ with $z=E(u)$ and $u_q=R(z)$, it holds that
\begin{equation}\label{eq:cdet-bound}
\|P_h(u)-\widehat P[W^\ast](u)\|_u\le C_L\,c_+\Big(\|\Pi\|_{\rm ref}\,\|P_h(u)-P_h(u_q)\|_{\rm ref}+|c(z)-W^\ast B(z)|\Big),
\end{equation}
and
\[
\big\|\,|c-W^\ast B|\,\big\|_{L^2(\Xi)}\le C\,M^{-\frac12-\frac{2k_\br+1}{2q}}\,\|c\|_{H^{r_{\max}(q,k_\br)}(\Xi;\R^{n_V})}.
\]
If, in addition, $\|\cdot\|_u$ extends to a seminorm on $\mathcal G(u)+V_h$, then
\begin{align*}
\|\mathcal G(u)-\widehat P[W^\ast](u)\|_u\le\ &\|\mathcal G(u)-P_h(u)\|_u+C_L c_+\|\Pi\|_{\rm ref}\,\|P_h(u)-P_h(u_q)\|_{\rm ref}\\
&+C_L c_+\,|c(z)-W^\ast B(z)|.
\end{align*}
\end{theorem}

\begin{proof}
Let $P_{h,n}=\Pi P_h(u)$ and $P_n=\sum_j(W^\ast B(z))_j\phi_j$.
By \cref{lem:layer}, the uniform layer bound, and the triangle inequality, we have
\[
\|P_h(u)-\widehat P[W^\ast](u)\|_u\le C_L c_+\big(\|\Pi(P_h(u)-P_h(u_q))\|_{\rm ref}+\|\Pi P_h(u_q)-P_n\|_{\rm ref}\big).
\]
Since $\Pi P_h(u_q)-P_n=\sum_j\big(c_j(z)-(W^\ast B(z))_j\big)\phi_j$ and $\{\phi_j\}$ is orthonormal, the last norm equals $|c(z)-W^\ast B(z)|$.
The $L^2(\Xi)$ bound follows from \cref{lem:branch} with $H=\R^{n_V}$ and $s=r_{\max}(q,k_\br)$, and the final bound follows from the triangle inequality.
This completes the proof.
\end{proof}

\subsection{Training without solution data and an error certificate}\label{sec:training}
Since $\widehat P[W]$ is affine in $W$, both the \emph{data loss} $\frac1{n_s}\sum_m\|P_h(u_m)-\widehat P[W](u_m)\|_{u_m}^2$, which requires the discrete solutions $P_h(u_m)$, and the \emph{residual loss} $\frac1{n_s}\sum_m\|\rho_m\|_{*,u_m}^2$, which requires no solutions, lead to linear least squares problems for $W$, given training inputs $u_1,\dots,u_{n_s}$.
Here, $\rho_m(w):=F_{u_m}(w)-a_{u_m}(\widehat P[W](u_m),w)$ is the residual and $\|\rho\|_{*,u}:=\sup_{w\in V_h}\rho(w)/\|w\|_u$ is the dual norm.

\begin{proposition}[Residual certificate]\label{prop:certificate}
For every input $u$ and every readout $W$, let $\rho(w)=F_u(w)-a_u(\widehat P[W](u),w)$.
Then it holds that
\begin{equation}\label{eq:certificate}
C_a^{-1}\|\rho\|_{*,u}\le\|P_h(u)-\widehat P[W](u)\|_u\le\|\rho\|_{*,u}.
\end{equation}
Moreover, $\rho(\psi)=0$ for all $\psi\in V_0$.
\end{proposition}

\begin{proof}
Let $e=P_h(u)-\widehat P[W](u)$.
Then $a_u(e,w)=\rho(w)$ for all $w\in V_h$.
By~\eqref{eq:coerc}, we have $\|e\|_u^2\le a_u(e,e)=\rho(e)\le\|\rho\|_{*,u}\|e\|_u$ and $\rho(w)\le C_a\|e\|_u\|w\|_u$.
The last claim is~\eqref{eq:layer-gen}.
This completes the proof.
\end{proof}

Hence, the residual loss is equivalent to the data loss sample by sample, and $\|\rho\|_{*,u}$ is a computable error certificate for every new input, which costs one solve with the Gram matrix of $\|\cdot\|_u$ on $V_h$.
Since the certificate requires no solutions, a ridge parameter can also be selected by minimizing the mean certificate over validation inputs, so that the whole training procedure uses no solution data.

\subsection{A locally conservative DetONet for Darcy flow}\label{sec:cons}
We now apply the constrained DetONet to the ENN method of~\cite{kim2026locally} for~\eqref{eq:darcy}, so that every prediction is exactly locally conservative on every element, while training remains a single linear least squares problem.
Let $\mathcal T_h$ be a fixed, possibly coarse, auxiliary partition of $\Omega$ into elements $K$; throughout, only the dimension of the neural space is increased, and $h$ labels the partition, not a refinement parameter.
The ENN space is given by
\[
V_h:=\mathcal N_n^k+\mathcal M^0(\mathcal T_h),
\]
where
\[
\mathcal N_n^k:=\spn\{\sigma_k(\theta_j\cdot(x,1))\}_{j=1}^n,\qquad
\mathcal M^0(\mathcal T_h):=\spn\{\psi_K\}_{K\in\mathcal T_h},
\]
the neuron parameters $\theta_j$ are fixed and quasi-uniform on $\Sph^d$, and $\psi_K$ is the characteristic function of $K$.
The discrete operator $\mathcal G_h:\kappa\mapsto P_h(\kappa)$ is the ENN solution of~\cite{kim2026locally} in its nonsymmetric interior penalty form, i.e., $P_h(\kappa)\in V_h$ satisfies $A_1(P_h(\kappa),w)=F_1(w)$ for all $w\in V_h$ with
\begin{equation}\label{eq:A1}
\begin{aligned}
A_1(v,w):=\ &(\kappa\nabla v,\nabla w)_{\mathcal T_h}-\langle\avg{\kappa\nabla v},\jump{w}\rangle_{\mathcal E_h^1}+\langle\jump{v},\avg{\kappa\nabla w}\rangle_{\mathcal E_h^1}\\
&+\gamma\langle\hat h_e^{-1}\kappa_e\jump{v},\jump{w}\rangle_{\mathcal E_h^1}+\tau(\nabla\cdot\kappa\nabla v,\nabla\cdot\kappa\nabla w)_{\mathcal T_h}\\
&+\gamma_N\langle\hat h_e\,\kappa\nabla v\cdot n,\kappa\nabla w\cdot n\rangle_{\mathcal E_h^N},\\
F_1(w):=\ &(f,w-\tau\nabla\cdot\kappa\nabla w)_{\mathcal T_h}-\langle j_N,w+\gamma_N\hat h_e\,\kappa\nabla w\cdot n\rangle_{\mathcal E_h^N}\\
&+\langle p_D,\avg{\kappa\nabla w}\cdot n\rangle_{\mathcal E_h^D}+\gamma\langle\hat h_e^{-1}\kappa_e\,p_D,w\rangle_{\mathcal E_h^D}.
\end{aligned}
\end{equation}
Here, $\mathcal E_h^1$ is the set of interior and Dirichlet edges, $\mathcal E_h^D$ and $\mathcal E_h^N$ are the sets of Dirichlet and Neumann edges, $\avg{\cdot}$ and $\jump{\cdot}$ are the weighted average and the jump, $\kappa_e$ is the edge value of $\kappa$, and $\hat h_e:=\min(h_e,\delta)$ with the diameter $h_e$ of the edge $e$ and $\delta=n^{-1/d}$.
The penalty parameter $\gamma>0$, the Galerkin least-squares parameter $\tau>0$, and the Neumann least-squares parameter $\gamma_N\ge0$ are denoted by $\alpha$, $\tau$, and $\gamma$ in~\cite{kim2026locally}.
With the energy norm
\[
\begin{aligned}
\|v\|_{\EG,\kappa}^2:=\ &\|\kappa^{1/2}\nabla v\|_{\mathcal T_h}^2+\gamma\sum_{e\in\mathcal E_h^1}\hat h_e^{-1}\|\kappa_e^{1/2}\jump{v}\|_{0,e}^2\\
&+\tau\|\nabla\cdot\kappa\nabla v\|_{\mathcal T_h}^2+\gamma_N\sum_{e\in\mathcal E_h^N}\hat h_e\|\kappa\nabla v\cdot n\|_{0,e}^2,
\end{aligned}
\]
it holds that $A_1(v,v)=\|v\|_{\EG,\kappa}^2$ for all $v\in V_h$~\cite[Lemma~3.5]{kim2026locally}.
Hence, \eqref{eq:coerc} holds with $a_\kappa=A_1$, $F_\kappa=F_1$, and $\|\cdot\|_\kappa=\|\cdot\|_{\EG,\kappa}$.
We assume $k\ge2$ and, as in~\cite[(2.2)]{kim2026locally}, $\kappa\in W^{1,\infty}(\Omega)$ with $0<\kappa_-\le\kappa\le\kappa_+$, so that the $\tau$-terms are well defined on $V_h$.
For each fixed neural space, the continuity constant of $A_1$ involves the trace constant
\begin{equation}\label{eq:Ctr}
C_{\rm tr}(n,h):=\sup_{v\in\mathcal N_n^k,\,\nabla v\ne0}\frac{\sum_{e\in\mathcal E_h^1}\hat h_e\|\kappa^{1/2}\nabla v\|_{0,e}^2}{\|\kappa^{1/2}\nabla v\|_{\mathcal T_h}^2+\tau\|\nabla\cdot\kappa\nabla v\|_{\mathcal T_h}^2}
\end{equation}
of~\cite[Lemma~3.4]{kim2026locally}, namely $|A_1(v,w)|\le(1+2(C_{\rm tr}/\gamma)^{1/2})\|v\|_{\EG,\kappa}\|w\|_{\EG,\kappa}$~\cite[(3.30)]{kim2026locally}; $C_{\rm tr}(n,h)$ is finite for every $n$, but it is not known to be bounded uniformly in $n$, since $\mathcal N_n^k$ has no inverse inequality.
The discretization error $\|p(\kappa)-P_h(\kappa)\|_{\EG,\kappa}$ is bounded in~\cite[Theorems~3.9 and~3.10]{kim2026locally}; if $\Gamma_N\neq\emptyset$, the C\'ea-type estimate~\cite[Theorem~3.9]{kim2026locally} requires $\gamma_N>0$.
In the notation of \cref{sec:layer}, $V_n=\mathcal N_n^k$ is the trunk and $V_0=\mathcal M^0(\mathcal T_h)$ carries the constraints.

In block form with respect to $\mathcal N_n^k$ and $\mathcal M^0(\mathcal T_h)$, the ENN system reads
\[
\begin{pmatrix}\mathsf A_{nn}(\kappa)&\mathsf A_{n0}(\kappa)\\ \mathsf A_{0n}(\kappa)&\mathsf A_{00}(\kappa)\end{pmatrix}
\begin{pmatrix}\mathbf c\\ \boldsymbol\mu\end{pmatrix}=\begin{pmatrix}\mathsf F_n(\kappa)\\ \mathsf F_0(\kappa)\end{pmatrix}.
\]
The rows tested with $\psi_K$ are local mass balances: by~\cite[Lemma~4.1]{kim2026locally}, if $P\in V_h$ satisfies
\begin{equation}\label{eq:psiK}
A_1(P,\psi_K)=F_1(\psi_K)\qquad\forall\,K\in\mathcal T_h,
\end{equation}
then the ENN flux $\mathbf U$ of $P$, defined in~\cite[(4.1)]{kim2026locally}, satisfies $\int_{\partial K}\mathbf U\cdot n_K\,ds=\int_Kf\,dx$ on every element.
Note that the proof uses only~\eqref{eq:psiK}, so that conservation holds for any $P\in V_h$ that satisfies these rows, whatever its neural part is.
Therefore, the layer~\eqref{eq:layer-gen} is a \emph{conservation layer}: given predicted neural coefficients $\mathbf c$, it computes
\begin{equation}\label{eq:layer}
\boldsymbol\mu(\kappa)=\mathsf A_{00}(\kappa)^{-1}\big(\mathsf F_0(\kappa)-\mathsf A_{0n}(\kappa)\,\mathbf c\big),\qquad
\widehat P(\kappa)=\sum_jc_j\phi_j+\sum_K\mu_K(\kappa)\,\psi_K.
\end{equation}
Since $\nabla\psi_K=0$, the block $\mathsf A_{00}$ contains only the penalty terms, i.e., $(\mathsf A_{00})_{KK'}=\allowbreak\gamma\langle\hat h_e^{-1}\kappa_e\jump{\psi_{K'}},\jump{\psi_K}\rangle_{\mathcal E_h^1}$; it is a weighted graph Laplacian of the element adjacency, grounded at the Dirichlet edges, and symmetric positive definite when the mesh graph is connected to a Dirichlet face.
Note that $\mathsf A_{0n}$, $\mathsf A_{00}$, and $\mathsf F_0$ involve $\kappa$ only through edge integrals, and that the layer uses the true input $\kappa$, not its encoding.

\begin{proposition}[Exact local conservation]\label{prop:cons}
For every readout $W$ and every input $\kappa$, the conservative DetONet
\[
\widehat P[W](\kappa)=\mathcal L_\kappa\Big(\sum_j\big(WB(E(\kappa))\big)_j\phi_j\Big)
\]
is locally conservative, i.e., its ENN flux satisfies $\int_{\partial K}\mathbf U\cdot n_K\,ds=\int_Kf\,dx$ for all $K\in\mathcal T_h$, up to the round-off of the $|\mathcal T_h|\times|\mathcal T_h|$ solve.
\end{proposition}

\begin{proof}
By construction, $\mathsf A_{0n}\mathbf c+\mathsf A_{00}\boldsymbol\mu=\mathsf F_0$, which is~\eqref{eq:psiK}.
Then the result follows from~\cite[Lemma~4.1]{kim2026locally}.
This completes the proof.
\end{proof}

For the encoder, we use the point values of $\log\kappa$ of \cref{sec:model}, i.e., $z=E(\kappa)=\mathbf E\,(\log\kappa)(\hat{\mathbf x})\in\R^q$ are the least squares coefficients of $\log\kappa$ in a basis $e_1,\dots,e_q$ from $m$ sensors, and $R(z)=\kappa_z$.
The branch has $M$ $\mathrm{ReLU}^{k_\br}$ neurons with parameters fixed and quasi-uniform on $\Sph^q$, and the neural part is predicted in a basis of $\mathcal N_n^k$ of dimension $n_V\le n$ that is orthonormal in the EG norm of a reference conductivity $\kappa_0$, which serves as $\|\cdot\|_{\rm ref}$.
The following theorem shows that, on a fixed auxiliary partition, the layer constant $C_{\rm layer}$ of \cref{lem:layer} is bounded independently of the neural dimension.

\begin{theorem}[Fixed-mesh stability of the conservation layer]\label{thm:fixed-mesh-layer}
Let $\mathcal T_h$ be a fixed conforming, shape-regular polygonal or polyhedral partition with a Dirichlet-anchored element adjacency graph, and assume that the admissible conductivities satisfy $\kappa\in W^{1,\infty}(\Omega)$ and $0<\kappa_-\le\kappa\le\kappa_+$ uniformly.
Let $\delta=n^{-1/d}$ and $\hat h_e=\min(h_e,\delta)$, suppose that $\gamma>0$ and, if $\Gamma_N\neq\emptyset$, $\gamma_N>0$ in~\eqref{eq:A1} are independent of $n$, and that $\tau\ge c_\tau\delta^2$.
Then
\[
\sup_{\kappa}\sup_{n\ge1} C_{\rm layer}(\kappa,n;\mathcal T_h)\le C_{\mathcal T_h}<\infty,
\]
where $C_{\mathcal T_h}$ may depend on the fixed partition, the uniform conductivity bounds, and the stabilization parameters, but not on the neural dimension $n$ and resolution $\delta$.
\end{theorem}

The proof, given in \cref{sm:sec-fixed}, bounds the skew coupling $D_\kappa$ of \cref{lem:layer} by elementwise Green's identity with cutoffs of width $O(\delta)$ along $\partial K$ and a graph Poincar\'e inequality on the fixed partition.
Note that \cref{thm:fixed-mesh-layer} bounds only this coupling, not the continuity constant of the full form, requires no inverse inequality in $n$, and asserts no limit $h\to0$.

In \cref{thm:cdet}, we take for $\Pi$ the projection onto $V_n$ along the $\|\cdot\|_{\rm ref}$-orthogonal complement of $V_n\cap V_0$ in $V_0$; this complement and $V_n$ form a direct sum equal to $V_h$.
Since the functions of $V_n$ are continuous, $V_n\cap V_0$ contains at most the constants.
Then $\|\Pi\|_{\rm ref}=(1-c_F^2)^{-1/2}$, where $c_F<1$ is the cosine of the Friedrichs angle between $V_n$ and $V_0$ in $\|\cdot\|_{\rm ref}$.
Note that \cref{thm:fixed-mesh-layer} does not bound $c_F$ uniformly in $n$; we therefore keep a uniform bound on $\|\Pi\|_{\rm ref}$ as an assumption, and examine it numerically in \cref{sec:num-cons}.
Then \cref{thm:cdet,prop:cons,thm:fixed-mesh-layer} give the following result.

\begin{corollary}[Conservative DetONet for Darcy flow]\label{cor:cons}
Assume that~\eqref{eq:kappa-bound} holds, including the reconstructions $\kappa_q$, that the EG norms of these conductivities are uniformly equivalent to the reference norm, that the fixed-mesh stabilization hypotheses of \cref{thm:fixed-mesh-layer} hold, and that $\|\Pi\|_{\rm ref}\le C_\Pi$ for all $n$.
Then, for every input $\kappa$ with $z=E(\kappa)$, the conservative DetONet with the readout $W^\ast$ of \cref{thm:cdet} is exactly locally conservative, and it holds that
\begin{multline*}
\|p(\kappa)-\widehat P[W^\ast](\kappa)\|_{\EG,\kappa}
\le\|p(\kappa)-P_h(\kappa)\|_{\EG,\kappa}\\
+C\,\|P_h(\kappa)-P_h(\kappa_q)\|_{\EG,\kappa_0}
+C\,|c(z)-W^\ast B(z)|,
\end{multline*}
with $C$ depending on the fixed auxiliary mesh, the uniform conductivity and stabilization bounds, the norm equivalence constants, and $C_\Pi$, but not on the neural dimension $n$.
The three terms are the ENN discretization error, the modeling error of $E$, and the branch error, which is of order $M^{-\frac12-\frac{2k_\br+1}{2q}}$ in $L^2(\Xi)$ by \cref{thm:cdet} if $c\in H^{r_{\max}(q,k_\br)}(\Xi;\R^{n_V})$.
Note that $\|c\|_{H^{r_{\max}(q,k_\br)}(\Xi;\R^{n_V})}$ may depend on $n$, through the dependence of $P_h$ on $\kappa$ (see \cref{lem:perturb}).
\end{corollary}

The uniform equivalence of the EG norms in \cref{cor:cons} holds, for example, if $\|\nabla\log\kappa\|_{L^\infty(\Omega)}$ is also uniformly bounded and $\tau\le1$, since $\nabla\cdot\kappa\nabla v=\kappa(\Delta v+\nabla\log\kappa\cdot\nabla v)$ on each element.
The modeling term $\|P_h(\kappa)-P_h(\kappa_q)\|$ in \cref{cor:cons} is the discrete counterpart of \cref{lem:darcy}.
For admissible $\kappa,\kappa'$, let $\omega:=1-\kappa'/\kappa$ and
\[
\eta:=\|\omega\|_{L^\infty(\Omega)}\le e^{\varepsilon}-1,\qquad
\eta_1:=\tau^{1/2}\|\kappa\|_{L^\infty(\Omega)}^{1/2}\|\nabla\omega\|_{L^\infty(\Omega)},
\]
where $\varepsilon=\|\log\kappa-\log\kappa'\|_{L^\infty(\Omega)}$, and write $P_h(\kappa,f)$ for the ENN solution with the source $f$.

\begin{lemma}[Coefficient perturbation of the ENN solution]\label{lem:perturb}
Let $\kappa,\kappa'$ be admissible, and let $C_{\rm tr}=C_{\rm tr}(n,h)$ be the trace constant~\eqref{eq:Ctr} for $\kappa$.
Then it holds that
\[
\|P_h(\kappa,f)-P_h(\kappa',f)\|_{\EG,\kappa}\le L(\eta,\eta_1)\,|||P_h(\kappa',f)|||_\kappa
\]
with
\[
L(\eta,\eta_1):=\eta\Big(1+\gamma^{1/2}+\big(1+\gamma^{-1/2}\big)C_{\rm tr}^{1/2}+\gamma_N^{1/2}(2+\eta)\Big)+(\eta+\eta_1)(2+\eta+\eta_1),
\]
where $|||\cdot|||_\kappa$ is the EG-type energy~\eqref{sm:triple} that includes the data.
\end{lemma}

The proof is given in \cref{sm:sec-perturb}, together with the corresponding estimate for an encoded source in a discrete dual norm.
As in the continuous bound~\eqref{eq:darcy-stab}, $L(\eta,\eta_1)$ is linear in $\eta$ and $\eta_1$ for small perturbations.
However, $L$ depends on $n$ through $C_{\rm tr}(n,h)$, which is not known to be bounded uniformly in $n$, and it contains $\eta_1$, since the $\tau$-term of~\eqref{eq:A1} involves $\nabla\kappa$; note that $\eta_1$ carries the factor $\tau^{1/2}$, which is of order $\delta$ for $\tau\simeq\delta^2$.
With $\kappa'=\kappa_q$, $\eta$ is bounded by \cref{lem:enc}, whereas $\eta_1$ requires an approximation of $\nabla\log\kappa$ by the encoder, which \cref{lem:enc} does not provide.

Finally, the readout is computed by one linear least squares solve with the data loss or the residual loss of \cref{sec:training}, and a prediction requires one encoding, one branch evaluation, the product $WB(z)$, the assembly of $\mathsf A_{0n}$, $\mathsf A_{00}$, and $\mathsf F_0$ from edge integrals, and one $|\mathcal T_h|\times|\mathcal T_h|$ solve; no PDE is solved.

\begin{remark}[Local versus strong conservation]\label{rem:strong}
The layer enforces conservation in the cellwise sense of~\eqref{eq:psiK}, whereas some neural approximations enforce it in the strong sense; for example, the divergence-free linearized neural networks of He, Liu, and Tian~\cite{he2026divergencefree} use antisymmetric potentials to enforce incompressibility exactly.
A potential-based flux $\mathbf U_N=\mathbf U_f+\nabla\times\mathbf A_N$ with $\nabla\cdot\mathbf U_f=f$ satisfies $\nabla\cdot\mathbf U_N=f$ in the sense of distributions, but it need not satisfy the constitutive law $\mathbf U_N=-\kappa\nabla q$ for some $q$ or the Neumann condition $\mathbf U_N\cdot n=j_N$ on $\Gamma_N$, so that strong conservation alone does not imply convergence to the solution of~\eqref{eq:darcy}.
The ENN layer instead obtains exact cellwise balance from the variational formulation~\eqref{eq:A1}, which also imposes the constitutive law and the boundary conditions weakly, so that conservation and the error estimate of \cref{cor:cons} hold simultaneously.
\end{remark}

\section{Numerical experiments}\label{sec:numerics}

We test network convergence with exact encoding, quantify conductivity and source encoding errors, and verify the locally conservative DetONet.

In all experiments, we consider the Darcy flow~\eqref{eq:darcy} on $\Omega=(0,1)^2$ with coordinates $x=(x_1,x_2)$, $p=1$ on $x_1=0$, $p=0$ on $x_1=1$, and no flux on $x_2\in\{0,1\}$.
The trunk and branch networks use $\mathrm{ReLU}^3$ neurons, i.e., $k_\tr=k_\br=3$, with quasi-uniform parameters: a Fibonacci set on $\Sph^2$, equispaced points on $\Sph^1$, and a set obtained by minimizing the Riesz energy on $\Sph^q$ for $q\ge3$.
As in~\cite{kim2026locally}, the trunk neurons are centered at $x_c=(\frac12,\frac12)$, i.e., $T_j(x)=\sigma_3(\theta^\tr_j\cdot(x-x_c,1))$, which is the DetONet of \cref{sec:detonet-def} on the shifted domain $\Omega-x_c$.
For comparison, we also use i.i.d.\ random parameters uniformly distributed on the same spheres.
All experiments are implemented in Python using NumPy and SciPy, and are conducted on a workstation equipped with an Apple M5 Max processor (18 cores) and 64\,GB of memory.

\subsection{Network error}\label{sec:num-rates}
First, we verify the rates of \cref{thm:rate}.
The conductivity is given by
\[
\log\kappa_\xi(x)=0.4\sum_{\ell=1}^qc_\ell(\xi)\cos(\pi a_\ell x_1)\cos(\pi b_\ell x_2),\qquad \xi\sim U[-1,1]^q,
\]
with $(a_\ell,b_\ell)$ taken in the order $(1,1),(2,1),(1,2),(2,2)$, and $f=0$.
We consider an \emph{analytic dependence}, $c_\ell(\xi)=\xi_\ell$ with $q\in\{2,4\}$, and a dependence of \emph{limited smoothness}, $q=2$, $c_1(\xi)=|\xi_1|^\alpha$, and $c_2(\xi)=\xi_2$ with $\alpha\in\{0.5,1.5\}$, so that $g\in H^s(\Xi;L^2(\Omega))$ exactly for $s<\alpha+\frac12$.
Since $\kappa$ is even about the sides of $\Omega$ and $f=0$, odd reflections of $p-1$ at $x_1=0$ and of $p$ at $x_1=1$ and even reflections at $x_2\in\{0,1\}$ are again solutions, so that $p$ is smooth up to the boundary; in particular, $p(z,\cdot)\in H^{r_{\max}(2,3)}(\Omega)=H^{4.5}(\Omega)$, and the trunk can reach its saturation rate.
Since $c_1(\xi)=|\xi_1|^\alpha$ does not distinguish $\pm\xi_1$, we regard the parameter $\xi$ as the input, i.e., $\mathcal U=[-1,1]^q$ and $\mathcal G(\xi)=p(\kappa_\xi)$, and the encoder returns $z=\xi/\sqrt q$; then $E$ is a bijection onto $\Xi=[-1/\sqrt q,1/\sqrt q]^q$, so that $\mathcal G(u)=g(z)$ and $\varepsilon_E=0$.

Reference solutions use a $128\times128$ cell-centered finite volume grid, so the target is the discrete operator. We use $4000$ training inputs for analytic $q=2$, $16000$ for $q=4$, $12000$ for limited smoothness, and $200$ test inputs. The readout uses an orthonormal trunk basis and branch SVD (cutoff $10^{-12}$); for $t=1$, the trunk readout uses the discrete $H^1$ inner product. Trunk tests fix a sufficiently rich branch ($M=512$ for $q=2$, $256$ for $q=4$), while branch tests represent the output exactly in space. Errors are relative $L^2$ or discrete $H^1$ errors over test inputs.

\begin{figure}[t]
\centering
\includegraphics[width=\textwidth]{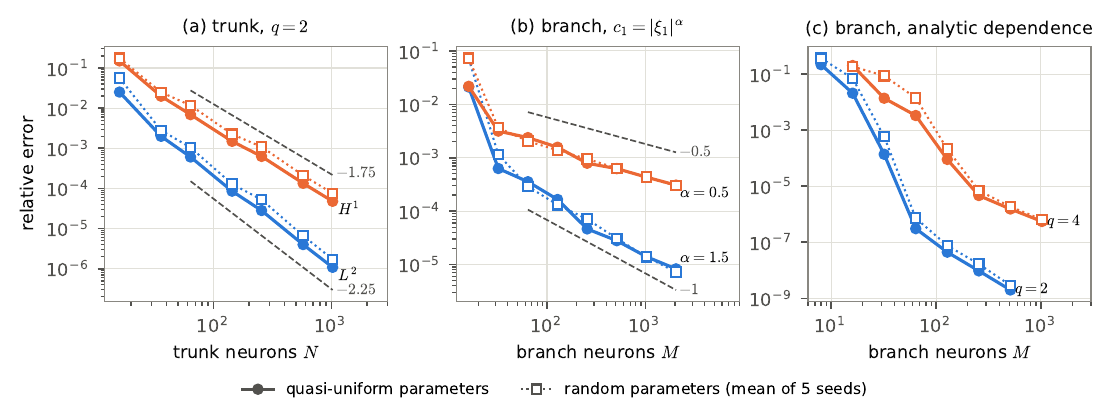}
\caption{Convergence of the network error with $\varepsilon_E=0$. (a) Trunk: relative $L^2$ and $H^1$ errors versus the number $N$ of trunk neurons for $q=2$. (b) Branch for the inputs of limited smoothness $c_1=|\xi_1|^\alpha$: relative $L^2$ error versus the number $M$ of branch neurons. (c) Branch for analytic dependence on $z$, $q\in\{2,4\}$. Solid lines with filled markers: quasi-uniform parameters; dotted lines with hollow markers: i.i.d.\ random parameters, mean over $5$ seeds. Dashed lines: the slopes predicted by \cref{thm:rate}.}\label{fig:rates}
\end{figure}

The observed trunk $L^2$ rates are $2.27$ ($q=2$) and $2.19$ ($q=4$), versus the predicted saturation exponent $2.25$; the $H^1$ rate is $1.78$, versus $1.75$ (\cref{fig:rates}). All rates are least squares fits over $N\in[64,1024]$ ($[64,576]$ for $q=4$) and $M\in[64,2048]$, and are listed in \cref{sm:tab-rates}. The $q=4$ trunk error eventually reaches the branch floor. For limited-smoothness branch targets, fitted rates $1.10$ ($\alpha=1.5$) and $0.59$ ($\alpha=0.5$) agree with predicted rates $1$ and $0.5$. Analytic targets converge faster than the saturation exponent over the tested range, which does not contradict~\cite{mao2025sharp}, since the saturation exponent is a worst-case rate over Sobolev balls. Random features give similar rates but trunk errors $1.5$--$1.9$ times larger and $17\%$--$59\%$ variability; quasi-uniform features provide deterministic guarantees. The discrete $H^1$ projection is the one directly covered by \cref{lem:trunk}; an $L^2$ readout has similar observed $H^1$ rates but is not covered without an inverse estimate.

\subsection{Modeling error}\label{sec:num-model}
Next, we examine the modeling error of \cref{sec:model} for an encoder that loses information, first for the conductivity
\[
\log\kappa(x)=\sum_{\ell=1}^{16}\xi_\ell\,c_\ell\cos(\pi a_\ell x_1)\cos(\pi b_\ell x_2),\qquad c_\ell=\frac{0.8}{a_\ell^2+b_\ell^2},\qquad\xi\sim U[-1,1]^{16},
\]
where $(a_\ell,b_\ell)\ne(0,0)$ are the $16$ wave numbers with $0\le a_\ell,b_\ell\le4$ of smallest $a_\ell^2+b_\ell^2$, ordered by $a_\ell^2+b_\ell^2$ and then lexicographically, and $f=0$.
The encoder is the least squares fit $\mathbf E=\mathbf S^\dagger$ of the first $q$ modes to $16\times16$ sensors of $\log\kappa$, so that \cref{lem:enc} applies, and the remaining $16-q$ modes are lost.
For each test input, we compute $\varepsilon_{\log}=\|\log\kappa-\log\kappa_q\|_{L^\infty(\Omega)}$ on a $201\times201$ grid and the solution $p(\kappa_q)$ for the truncated conductivity.
The DetONet is exact in $x$ as in \cref{sec:num-rates}, its branch has $M$ neurons acting on $z/\bar z$, where $\bar z$ is $1.05$ times the largest $|z|$ in the training set, and it is trained on $n_s=8000$ inputs.
The $L^2$ errors are relative errors over all test inputs, i.e., $(\sum_m\|e_m\|_{L^2}^2/\sum_m\|p_m\|_{L^2}^2)^{1/2}$ for the errors $e_m$ and the solutions $p_m$, whereas the other quantities are means over the test inputs.
We add a ridge term and select $\lambda\in\{0,10^{-12},10^{-11},\dots,10^{-2}\}$ by fitting on $7000$ of the training inputs and validating on the remaining $1000$, after which the readout is recomputed on all $n_s$ inputs.

\begin{figure}[t]
\centering
\includegraphics[width=\textwidth]{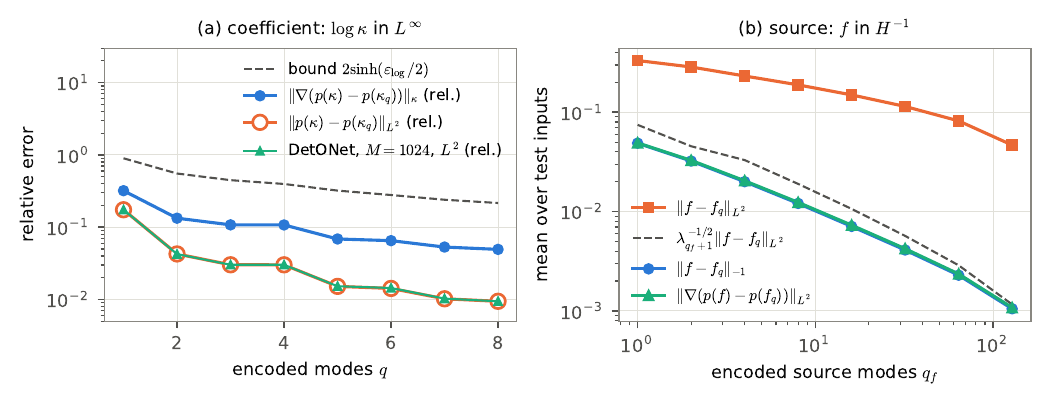}
\caption{Modeling error of the encoder. (a) Conductivity: the relative energy error $\|\kappa^{1/2}\nabla(p(\kappa)-p(\kappa_q))\|/\|\kappa_q^{1/2}\nabla p(\kappa_q)\|$ with the bound $2\sinh(\varepsilon_{\log}/2)$ of~\eqref{eq:darcy-sinh}, and the relative $L^2$ errors of $p(\kappa_q)$ and of the DetONet with $M=1024$, versus the number $q$ of encoded modes. (b) Source: the $L^2$ and $H^{-1}_D$ norms of $f-f_q$ and the error $\|\nabla(p(f)-p(f_q))\|_{L^2}$ (means over $100$ test inputs) versus the number $q_f$ of encoded source modes; dashed: the bound~\eqref{eq:negnorm}.}\label{fig:model}
\end{figure}

\begin{table}[t]
\centering
\footnotesize
\setlength{\tabcolsep}{4pt}
\caption{Modeling error of the encoder for $\log\kappa$ ($200$ test inputs): the Lebesgue constant $\Lambda$, the mean encoding error $\varepsilon_{\log}$, the mean relative energy error of $p(\kappa_q)$ with the mean of the bound $2\sinh(\varepsilon_{\log}/2)$, and the relative $L^2$ errors of $p(\kappa_q)$ and of the DetONet with the validated ridge parameter and with $\lambda=0$.}\label{tab:model}
\begin{tabular}{rcccccccc}
\toprule
 & & & \multicolumn{3}{c}{Modeling error} & \multicolumn{3}{c}{DetONet, relative $L^2$} \\
\cmidrule(lr){4-6}\cmidrule(lr){7-9}
$q$ & $\Lambda$ & $\varepsilon_{\log}$ & energy & $2\sinh(\varepsilon_{\log}/2)$ & $L^2$ & $M=256$ & $M=1024$ & $M=1024$, $\lambda=0$ \\
\midrule
1 & 1.28 & 0.861 & 3.20e-1 & 8.97e-1 & 1.74e-1 & 1.74e-1 & 1.74e-1 & 1.76e-1 \\
2 & 1.62 & 0.544 & 1.33e-1 & 5.52e-1 & 4.26e-2 & 4.25e-2 & 4.24e-2 & 4.46e-2 \\
3 & 2.20 & 0.443 & 1.08e-1 & 4.47e-1 & 3.02e-2 & 3.02e-2 & 3.02e-2 & 3.17e-2 \\
4 & 2.43 & 0.393 & 1.07e-1 & 3.96e-1 & 3.01e-2 & 3.01e-2 & 3.02e-2 & 3.26e-2 \\
5 & 2.57 & 0.319 & 6.86e-2 & 3.21e-1 & 1.52e-2 & 1.52e-2 & 1.52e-2 & 2.75e-2 \\
6 & 2.83 & 0.277 & 6.51e-2 & 2.79e-1 & 1.42e-2 & 1.44e-2 & 1.44e-2 & 1.66e-2 \\
7 & 2.91 & 0.239 & 5.31e-2 & 2.40e-1 & 1.02e-2 & 1.04e-2 & 1.03e-2 & 1.13e-2 \\
8 & 2.94 & 0.215 & 4.93e-2 & 2.16e-1 & 9.42e-3 & 1.10e-2 & 9.43e-3 & 1.02e-2 \\
\bottomrule
\end{tabular}

\end{table}

The energy error stays below $2\sinh(\varepsilon_{\log}/2)$, as predicted by~\eqref{eq:darcy-sinh} (\cref{fig:model,tab:model}). The encoder Lebesgue constant ranges from $1.28$ to $2.94$. With $M=1024$ and ridge regularization, the DetONet error is within $1.5\%$ of the $L^2$ modeling error for every $q$, and for $q=2$ it is even slightly smaller ($4.24\cdot10^{-2}$ versus $4.26\cdot10^{-2}$). This is expected, since the least squares readout approximates the conditional mean of $p(\kappa)$ given $z$, whose mean-square error can be smaller than that of $p(\kappa_q)$. Without regularization, the error can be $1.8$ times larger. Thus, empirical least squares may require stabilization even when the continuous approximation bound is sharp.

Next, in order to show that the source is truncated in the negative norm, let $\log\kappa(x)=0.4\big(\xi_1\cos(\pi x_1)\cos(\pi x_2)+\xi_2\cos(2\pi x_1)\cos(\pi x_2)\big)$ with $\xi\sim U[-1,1]^2$, and
\[
f(x)=\sum_{j=1}^{200}\nu_j\,c_j\,\varphi_j(x),\qquad \varphi_j(x)=\sin(\pi a_jx_1)\cos(\pi b_jx_2),\qquad c_j=(\lambda_j/\lambda_1)^{-3/4},
\]
with $\nu\sim U[-1,1]^{200}$, where $\varphi_j$ are the eigenfunctions of $-\Delta$ with the boundary conditions of $H^1_D(\Omega)$, ordered by their eigenvalues $\lambda_j=\pi^2(a_j^2+b_j^2)$.
The encoder keeps the first $q_f$ coefficients of the source, so that $f_q$ is the $L^2$ projection of \cref{sec:model}.
Note that $\beta=0.8$ in~\eqref{eq:kappa-bound}, so that \cref{lem:darcy} gives $\|\nabla(p(f)-p(f_q))\|_{L^2}\le e^{0.8}\|f-f_q\|_{-1}\approx2.23\,\|f-f_q\|_{-1}$.

In \cref{fig:model}(b), the solution error follows $\|f-f_q\|_{-1}$: the ratio of the error to $\|f-f_q\|_{-1}$ is at most $1.17$ for every test input, and the fitted rate is $0.83$ over $q_f\in[4,128]$ versus $0.44$ for $\|f-f_q\|_{L^2}$. Hence, the bound $e^\beta\|f-f_q\|_{-1}$ of~\eqref{eq:darcy-model} overestimates the mean error by a factor of about $2.2$ for every $q_f$, and every individual error by at least $e^{0.8}/1.17\approx1.9$, whereas the bound $e^\beta\lambda_1^{-1/2}\|f-f_q\|_{L^2}$ obtained from the $L^2$ norm and the Poincar\'e inequality overestimates the mean error by factors from $4.8$ at $q_f=1$ to $31$ at $q_f=128$. This confirms~\eqref{eq:negnorm}.

\subsection{Locally conservative DetONet}\label{sec:num-cons}
Finally, we test the locally conservative DetONet of \cref{sec:cons}.
The conductivity is given by
\[
\log\kappa(x)=\sum_{k\in\mathcal K}a_k\big(\xi_k\cos(2\pi k\cdot x)+\zeta_k\sin(2\pi k\cdot x)\big),\qquad \xi_k,\zeta_k\sim N(0,1),
\]
where $\mathcal K$ is the set of the $24$ wave vectors $k\in\mathbb Z^2$ with $|k_1|,|k_2|\le3$ and either $k_1>0$ or $k_1=0<k_2$, so that $\log\kappa$ has $48$ modes, the amplitudes $a_k\propto(1+|k|^2)^{-1}$ are normalized to a pointwise standard deviation of $0.5$, and $f=0$.
Since the coefficients are Gaussian, $\log\kappa$ is not uniformly bounded; the analysis of \cref{sec:cons} applies to the sampled inputs, which satisfy~\eqref{eq:kappa-bound} with $\beta=2.74$ together with their reconstructions $\kappa_q$, and all test codes lie within the range of the training codes.
For the ENN space, we follow the setting of~\cite[Section~5.1]{kim2026locally}: $n=144$ $\mathrm{ReLU}^3$ neurons with a Fibonacci set on $\Sph^2$ and the center $x_c=(\frac12,\frac12)$, an $8\times8$ grid of square elements, the NIPG form~\eqref{eq:A1} with $\gamma=10$ and $\tau=\delta^2$, where $\delta=n^{-1/2}=1/12$, so that $\hat h_e=\delta$, and a composite three-point Gauss rule on a uniform background grid of $96\times96$ cells for all integrals.
Unlike the experiments of~\cite{kim2026locally}, which use $\gamma_N=0$, we use $\gamma_N=1$, so that the hypotheses of \cref{thm:fixed-mesh-layer} and of~\cite[Theorems~3.9 and~3.10]{kim2026locally} hold.
Since $35$ of the $144$ neurons vanish on $\Omega$ and the remaining features are numerically dependent, the trunk $\{\phi_j\}$ is an orthonormal basis of $V_n$ in the EG norm of $\kappa_0=1$, computed from the eigenvalue decomposition of the Gram matrix with relative cutoff $10^{-12}$; its dimension is $n_V=84$.

The encoder is the least squares fit of \cref{lem:enc} to the $q\in\{2,4,8,16\}$ modes of largest amplitude from $16\times16$ sensors of $\log\kappa$, where modes of equal amplitude are ordered with the cosines before the sines and the wave vectors in lexicographic order, so that the remaining $48-q$ modes are lost.
The branch has $M$ $\mathrm{ReLU}^3$ neurons with parameters quasi-uniform on $\Sph^q$ acting on $z/\bar z$, where $\bar z$ is $1.05$ times the largest $|z|$ in the training set, and the branch features are whitened on the training inputs with relative eigenvalue cutoff $10^{-10}$.
The readout is computed by the conjugate gradient method on the normal equations of the residual loss or the data loss of \cref{sec:training}.
We use $n_s=16000$ training inputs, and the ridge parameter $\lambda\in\{0,10^{-4},10^{-3},3\cdot10^{-3},\dots,0.3\}$ is selected by the mean certificate $\|\rho\|_{*,\kappa}^2$ over $1000$ validation inputs, so that training with the residual loss uses no ENN solution at all; for comparison, we also report $n_s=4000$ without regularization.
The errors are the relative EG and $L^2$ errors against $P_h(\kappa)$, averaged over $200$ test inputs, the modeling error of $E$ is $\|P_h(\kappa)-P_h(\kappa_q)\|_{\EG,\kappa}/\|P_h(\kappa)\|_{\EG,\kappa}$, and the mass defect of a prediction $\widehat P$ is $\max_K|A_1(\widehat P,\psi_K)-F_1(\psi_K)|$.
The continuity constant of $A_1$ in~\eqref{eq:coerc} is $C_a=1.070$, maximized over the test inputs.

\begin{table}[t]
\centering
\small
\caption{Conservative DetONet: mean relative EG error (\%) against $P_h(\kappa)$ for encoders keeping $q$ of the $48$ modes of $\log\kappa$ and $M$ branch neurons.
Top: residual loss (no solutions), $n_s=16000$, $\lambda$ selected by the certificate.
Middle: data loss with the same $\lambda$, residual loss with $n_s=4000$ and $\lambda=0$, and the modeling error of $E$.
Bottom: mean encoding error $\varepsilon_{\log}$ and Lebesgue constant $\Lambda$.}\label{tab:cdet}
\begin{tabular}{lcccc}
\toprule
 & $q=2$ & $q=4$ & $q=8$ & $q=16$ \\
\midrule
$M=32$ & 3.93 & 3.43 & 3.86 & 4.55 \\
$M=64$ & 3.93 & 3.39 & 3.35 & 3.98 \\
$M=128$ & 3.93 & 3.40 & 2.82 & 2.74 \\
$M=256$ & 3.94 & 3.41 & 2.76 & 2.30 \\
$M=512$ & 3.94 & 3.46 & 2.79 & 2.15 \\
$M=1024$ & 3.94 & 3.51 & 2.87 & 2.11 \\
\midrule
data loss, $M=1024$ & 3.94 & 3.51 & 2.88 & 2.11 \\
$\lambda=0$, $n_s=4000$, $M=1024$ & 3.97 & 4.02 & 3.57 & 2.61 \\
modeling error of $E$ & 3.92 & 3.37 & 2.72 & 2.01 \\
\midrule
$\varepsilon_{\log}$ & 1.15 & 0.98 & 0.77 & 0.57 \\
$\Lambda$ & 1.61 & 1.63 & 2.21 & 2.88 \\
\bottomrule
\end{tabular}

\end{table}

\begin{table}[t]
\centering
\footnotesize
\setlength{\tabcolsep}{3pt}
\caption{Effect of the conservation layer for $q=16$ and $M=256$, $n_s=16000$. ``Solutions'' indicates whether ENN solutions are used in training. The network without the layer predicts all ENN coefficients from the same branch by least squares in the EG norm of $\kappa_0$; the second row applies the layer to its neural part after training. The last row replaces the branch by $400$ random Fourier features of the $256$ sensor values, as in RandONets~\cite{fabiani2025randonets}. Effectivity: $\|\rho\|_{*,\kappa}$ over the EG error, which lies in $[1,C_a]$ by~\eqref{eq:certificate} for every prediction in $V_h$.}\label{tab:layer}
\resizebox{\textwidth}{!}{\begin{tabular}{lccccc}
\toprule
Model & Solutions & EG (\%) & $L^2$ (\%) & Mass defect (mean / max) & Effectivity \\
\midrule
DetONet without layer (data, coefficients) & yes & 3.19 & 2.84 & $6.3\cdot10^{-1}$ / $6.3\cdot10^{0}$ & [1.007, 1.045] \\
same readout, layer after training & yes & 2.31 & 1.80 & $2.0\cdot10^{-13}$ / $1.6\cdot10^{-12}$ & [1.004, 1.016] \\
conservative DetONet, data loss & yes & 2.30 & 1.83 & $4.1\cdot10^{-13}$ / $2.5\cdot10^{-12}$ & [1.004, 1.017] \\
conservative DetONet, residual loss & no & 2.30 & 1.83 & $4.4\cdot10^{-13}$ / $2.3\cdot10^{-12}$ & [1.004, 1.017] \\
RFF branch on sensors, layer, residual loss & no & 1.19 & 1.56 & $8.9\cdot10^{-13}$ / $4.1\cdot10^{-12}$ & [1.006, 1.024] \\
\bottomrule
\end{tabular}
}
\end{table}

We first examine the constants of \cref{thm:fixed-mesh-layer,cor:cons} on the fixed $8\times8$ partition for $n\in\{36,\dots,1024\}$, maximized over $\kappa_0=1$ and $20$ random conductivities (\cref{sm:tab-Clayer}).
As predicted by \cref{thm:fixed-mesh-layer}, $C_{\rm layer}$ stays between $1.055$ and $1.065$ while $\delta$ decreases from $1/6$ to $1/32$, and it stays below $C_a$ for every $n$, in agreement with~\eqref{eq:layer-Ca}; without the Neumann least-squares term ($\gamma_N=0$), which is not covered by \cref{thm:fixed-mesh-layer}, it is almost the same.
Moreover, $\|\Pi\|_{\rm ref}$ increases only from $2.04$ to $2.10$, which supports the assumption of \cref{cor:cons}, and the trace constant $C_{\rm tr}$ of \cref{lem:perturb} decreases from $2.63$ to $2.02$ for $n\ge64$, in agreement with the values between $1.4$ and $2.7$ reported in~\cite[Section~4]{kim2026locally}.
For $\kappa'=\kappa_q$ and for smaller perturbations, the ratio of the error to the bound of \cref{lem:perturb} is at most $0.034$ (\cref{sm:sec-numerics}).

The conservation layer enforces cellwise mass balance to round-off: the largest defect is $4.2\cdot10^{-12}$ over all trained conservative DetONets and $2.0\cdot10^{-11}$ over every output of the layer, compared with the mean defect $6.3\cdot10^{-1}$ and the largest defect $6.3$ without the layer (\cref{tab:layer}).
The ratio $\|P_h(\kappa)-\widehat P(\kappa)\|_{\EG,\kappa}/\|P_{h,n}-P_n\|_{\EG,\kappa}$ of \cref{lem:layer}, with $P_{h,n}$ the neural part of the ENN solution, is at most $0.32$ for the trained networks and at most $1.01$ when the layer is applied after training, consistent with $C_{\rm layer}\approx1.06$.
Applying the layer after unconstrained training lowers the EG error from $3.19\%$ to $2.31\%$.

Residual training, requiring no ENN solutions, agrees with data training up to a relative difference of $0.02\%$, so that the values in \cref{tab:cdet} differ at most by rounding.
The computable certificate has effectivity between $1.003$ and $1.024$ for the trained networks and at most $1.045$ for every network, consistent with~\eqref{eq:certificate} and $C_a=1.070$.
As \cref{tab:cdet} and \cref{sm:fig-cdet} show, the error approaches the reconstruction-based modeling error of $E$ as $M$ increases: the best errors are $3.93\%$, $3.39\%$, $2.76\%$, and $2.11\%$ for $q=2,4,8,16$, versus the modeling errors $3.92\%$, $3.37\%$, $2.72\%$, and $2.01\%$.
For $q=4$ and $q=8$, it increases slightly for large $M$, by about $0.1$ percentage points at $M=1024$.
The selected ridge parameter was at most $10^{-4}$ and changed the errors by less than $10^{-3}$ percentage points, so that the improvement over $n_s=4000$ without regularization ($3.97\%$, $4.02\%$, $3.57\%$, and $2.61\%$ at $M=1024$) is due to the larger number of training inputs, which the residual loss provides without additional solves.
The table also shows the trade-off of \cref{cor:cons} between the modeling and branch errors: for $M=32$, the error increases with $q\ge4$ ($3.43\%$, $3.86\%$, and $4.55\%$), which is consistent with the branch exponent $\frac12+\frac{2k_\br+1}{2q}$ decreasing with $q$, although the constants also depend on $q$, whereas the best error decreases with $q$, as the modeling error does.

Finally, a random Fourier feature branch using all $256$ sensor values attains $1.19\%$ with the same conservation layer (\cref{tab:layer}), below the compressed DetONet.
It uses more input information than the $q$ encoded modes, so that it is not limited by the modeling error of $E$, but it has no corresponding approximation-rate guarantee.

\section{Conclusion}\label{sec:conclusion}

DetONet combines fixed-feature least-squares operator learning with exact elementwise conservation through an ENN correction layer. Continuous projection yields trunk and branch approximation rates; on a fixed auxiliary mesh the layer is stable independently of neural resolution under the stated hypotheses, including mixed boundary conditions. The residual supplies a discrete error certificate and permits solution-free training. Experiments support the rates, conservation to round-off, and a layer constant that is independent of the neural dimension on a fixed mesh. The projection estimates do not by themselves establish empirical least-squares convergence. Moreover, the constant of the coefficient-perturbation estimate depends on the neural trace constant, and the angle between the neural and constraint spaces enters the error bound; both were observed to be bounded independently of the neural dimension, but uniform bounds remain to be proved. Finally, the corrected network uses the full conductivity in its layer, so the encoding floor applies directly to the encoder-only prediction, not necessarily to every corrected architecture.

\appendix
\crefalias{section}{appendix}

\section{Proof of \cref{lem:darcy}}\label{sm:sec-darcy}

\begin{proof}[Proof of \cref{lem:darcy}]
We write $p(\kappa,f)-p(\kappa',f')=[p(\kappa,f)-p(\kappa',f)]+[p(\kappa',f)-p(\kappa',f')]$.
Each bracket is a difference of solutions with the same boundary data, so that it lies in $H^1_D(\Omega)$ and the Neumann terms cancel in its weak form.
First, the function $w:=p(\kappa',f)-p(\kappa',f')$ satisfies $(\kappa'\nabla w,\nabla v)=\langle f-f',v\rangle$ for all $v\in H^1_D(\Omega)$.
Testing with $v=w$ gives $e^{-\beta}\|\nabla w\|^2\le\|f-f'\|_{-1}\|\nabla w\|$.
Next, let $p=p(\kappa,f)$ and $p'=p(\kappa',f)$.
Subtracting the weak forms gives $(\kappa\nabla(p-p'),\nabla v)=((\kappa'-\kappa)\nabla p',\nabla v)$ for all $v\in H^1_D(\Omega)$.
Testing with $v=p-p'$, we obtain
\[
e^{-\beta}\|\nabla(p-p')\|^2\le\|\kappa-\kappa'\|_{L^\infty}\,\|\nabla p'\|\,\|\nabla(p-p')\|.
\]
Here, $\|\nabla p'\|\le C_{\nabla p}$, and $|e^a-e^b|\le e^{\max(a,b)}|a-b|$ gives $\|\kappa-\kappa'\|_{L^\infty}\le e^\beta\varepsilon$, which proves~\eqref{eq:darcy-stab}.
For~\eqref{eq:darcy-sinh}, we test the same identity with $v=p-p'$ and write $|\kappa'-\kappa|=\kappa^{1/2}\kappa'^{1/2}|e^{a/2}-e^{-a/2}|$ with $a=\log\kappa'-\log\kappa$.
Since $|e^{a/2}-e^{-a/2}|=2\sinh(|a|/2)\le2\sinh(\varepsilon/2)$, the Cauchy--Schwarz inequality gives the bound.
This completes the proof.
\end{proof}

\section{Proof of \cref{thm:fixed-mesh-layer}}\label{sm:sec-fixed}

By \cref{lem:layer}, it suffices to bound the skew coupling $D_\kappa$ between $\mathcal N_n^k$ and $\mathcal M^0(\mathcal T_h)$ uniformly in $n$.

\begin{proof}[Proof of \cref{thm:fixed-mesh-layer}]
Let $v\in\mathcal N_n^k$ and $\psi=\sum_K\mu_K\psi_K\in\mathcal M^0(\mathcal T_h)$ with $\mu_K\in\R$.
Since $v$ is continuous and $\nabla\psi=0$, the skew coupling of~\eqref{eq:A1} to the elementwise constants is
\[
(D_\kappa v)(\psi)=-\sum_{e\in\mathcal E_h^1}\int_e\kappa\nabla v\cdot\jump{\psi}\,ds.
\]
The sign convention for jumps does not affect the estimate.
Fix first $\delta$ smaller than a mesh-dependent positive threshold.
For each element, construct an $H^1(K)$ cutoff $\chi_K$ equal to $\mu_K$ on its entire boundary and supported in a strip of width $O(\delta)$.
Standard scaled distance-to-boundary cutoffs give
\[
\sum_K\|\chi_K\|_{0,K}^2\le C_{\mathcal T_h}\delta\, b(\psi)^2,\qquad
\sum_K\|\nabla\chi_K\|_{0,K}^2\le C_{\mathcal T_h}\delta^{-1}b(\psi)^2
\]
with $b(\psi)^2:=\sum_K|\partial K|\mu_K^2$.
Elementwise Green's identity yields
\[
(D_\kappa v)(\psi)=-\sum_K\big\{(\nabla\cdot\kappa\nabla v,\chi_K)_K+(\kappa\nabla v,\nabla\chi_K)_K\big\}+\langle\kappa\nabla v\cdot n,\psi\rangle_{\Gamma_N}.
\]
The anchored, fixed-mesh graph Poincar\'e inequality gives
\[
b(\psi)^2\le C_{\mathcal T_h}\sum_{e\in\mathcal E_h^1}\|\jump{\psi}\|_{0,e}^2.
\]
Since $\hat h_e=\delta$ in this regime, the penalty energy on $V_0$ implies $b(\psi)\le C_{\mathcal T_h}\sqrt{\delta}\,\|\psi\|_{\EG,\kappa}$.
Consequently, using $\tau\ge c_\tau\delta^2$, the two volume terms are bounded by
\[
C_{\mathcal T_h}\big(\|\kappa^{1/2}\nabla v\|_{\mathcal T_h}+\delta\|\nabla\cdot\kappa\nabla v\|_{\mathcal T_h}\big)\|\psi\|_{\EG,\kappa}\le C_{\mathcal T_h}\|v\|_{\EG,\kappa}\|\psi\|_{\EG,\kappa}.
\]
For mixed boundaries, the Neumann least-squares term of the EG norm and the same graph inequality give
\[
\sqrt{\delta}\,\|\kappa\nabla v\cdot n\|_{0,\Gamma_N}\le\gamma_N^{-1/2}\|v\|_{\EG,\kappa},\qquad
\delta^{-1/2}\|\psi\|_{0,\Gamma_N}\le C_{\mathcal T_h}\|\psi\|_{\EG,\kappa},
\]
respectively.
Thus, the Neumann term obeys the same bound; it vanishes for pure Dirichlet conditions.
We have proved
\[
|(D_\kappa v)(\psi)|\le C_{\mathcal T_h}\|v\|_{\EG,\kappa}\|\psi\|_{\EG,\kappa}.
\]
For the remaining values of $n$, $\delta=n^{-1/d}$ exceeds a fixed positive threshold, so that there are only finitely many such integers $n$.
The same estimate follows with a fixed-width boundary cutoff and the fixed-mesh graph inequality; all constants remain uniform over the admissible conductivity family.
Taking the dual supremum over $\psi\in V_0$ and applying \cref{lem:layer} proves the assertion.
This completes the proof.
\end{proof}

\begin{remark}[Scope of the stability statement]
The fixed-mesh result concerns the restricted conservation-layer coupling $D_\kappa$, not the continuity constant of the full nonsymmetric ENN form.
It does not require an inverse inequality uniform in the neural dimension.
It also does not establish a constant uniform under auxiliary mesh refinement.
For mixed boundary conditions, the proof uses the positive Neumann least-squares stabilization.
\end{remark}

\begin{remark}[Pure Neumann problems]
Although we assume $|\Gamma_D|>0$ elsewhere in this paper, the layer also applies if $\Gamma_D=\emptyset$.
In this case, the EG norm is a norm only on the functions of zero mean, and the constants span the kernel of $\mathsf A_{00}$.
Since $A_1(P_n,1)=0$ by the structure of~\eqref{eq:A1} and $F_1(1)=0$ by the compatibility condition, \eqref{eq:psiK} is solvable, and the layer solves~\eqref{eq:layer} in the subspace of zero mean.
\end{remark}

\section{Coefficient perturbation of the ENN solution}\label{sm:sec-perturb}

We write $A_1^\kappa$ and $F_1^\kappa$ to indicate the dependence on the conductivity, and $P_h(\kappa,f)$ for the ENN solution with the data $f$, $j_N$, and $p_D$.
For admissible $\kappa,\kappa'$, let $\omega$, $\eta$, and $\eta_1$ be as in \cref{sec:cons}, so that $\kappa-\kappa'=\omega\kappa$ and $\nabla\omega=-(\kappa'/\kappa)\nabla(\log\kappa'-\log\kappa)$; that is, $\eta_1$ measures the gradient of the perturbation of $\log\kappa$ with the factor $\tau^{1/2}$.
For $P\in V_h$, let $\jump{P}_D:=\jump{P}$ on interior edges and $\jump{P}_D:=(P-p_D)n$ on Dirichlet edges, and let
\begin{equation}\label{sm:triple}
\begin{aligned}
|||P|||_\kappa^2:=\ &\|\kappa^{1/2}\nabla P\|_{\mathcal T_h}^2+\sum_{e\in\mathcal E_h^1}\hat h_e^{-1}\|\kappa_e^{1/2}\jump{P}_D\|_{0,e}^2\\
&+\tau\big(\|\nabla\cdot\kappa\nabla P\|_{\mathcal T_h}^2+\|\nabla\cdot\kappa\nabla P+f\|_{\mathcal T_h}^2\big)\\
&+\sum_{e\in\mathcal E_h^N}\hat h_e\big(\|\kappa\nabla P\cdot n\|_{0,e}^2+\|\kappa\nabla P\cdot n+j_N\|_{0,e}^2\big)
\end{aligned}
\end{equation}
be an energy that includes the data.

\begin{proof}[Proof of \cref{lem:perturb}]
Let $P':=P_h(\kappa',f)$.
We first show that
\begin{equation}\label{sm:residual}
|F_1^\kappa(w)-A_1^\kappa(P',w)|\le L(\eta,\eta_1)\,|||P'|||_\kappa\,\|w\|_{\EG,\kappa}\qquad\forall\,w\in V_h.
\end{equation}
Let $\mathcal R(w):=\big(A_1^\kappa(P',w)-F_1^\kappa(w)\big)-\big(A_1^{\kappa'}(P',w)-F_1^{\kappa'}(w)\big)$, which equals $A_1^\kappa(P',w)-F_1^\kappa(w)$, since $A_1^{\kappa'}(P',w)=F_1^{\kappa'}(w)$.
By~\eqref{eq:A1}, the Dirichlet data enter $A_1(P',w)-F_1(w)$ only through $\jump{P'}_D$, the data $f$ and $j_N$ enter the $\tau$- and $\gamma_N$-terms only through $\nabla\cdot\kappa\nabla P'+f$ and $\kappa\nabla P'\cdot n+j_N$, and $(f,w)_{\mathcal T_h}$ and $\langle j_N,w\rangle_{\mathcal E_h^N}$ do not depend on $\kappa$.
Since $k\ge2$, $\nabla P'$ and $\nabla w$ are continuous, and since $\kappa_e=\kappa|_e$, we obtain
\[
\begin{aligned}
\mathcal R(w)=\ &(\omega\kappa\nabla P',\nabla w)_{\mathcal T_h}-\langle\omega\kappa\nabla P',\jump{w}\rangle_{\mathcal E_h^1}+\langle\jump{P'}_D,\omega\kappa\nabla w\rangle_{\mathcal E_h^1}\\
&+\gamma\langle\hat h_e^{-1}\omega\kappa_e\jump{P'}_D,\jump{w}\rangle_{\mathcal E_h^1}+\mathcal R_\tau(w)+\mathcal R_N(w),
\end{aligned}
\]
where $\mathcal R_\tau$ and $\mathcal R_N$ are the differences of the $\tau$- and $\gamma_N$-terms.
Since $|\omega|\le\eta$, the Cauchy--Schwarz inequality bounds the first and the fourth terms by $\eta\,|||P'|||_\kappa\|w\|_{\EG,\kappa}$ and $\eta\gamma^{1/2}|||P'|||_\kappa\|w\|_{\EG,\kappa}$, respectively.
For the second and the third terms, we use the weights $\hat h_e^{\pm1/2}$ and~\eqref{eq:Ctr}, i.e., $\sum_{e\in\mathcal E_h^1}\hat h_e\|\kappa^{1/2}\nabla v\|_{0,e}^2\le C_{\rm tr}(\|\kappa^{1/2}\nabla v\|_{\mathcal T_h}^2+\tau\|\nabla\cdot\kappa\nabla v\|_{\mathcal T_h}^2)$ for $v\in V_h$, which holds since the piecewise constants do not contribute to $\nabla v$, together with $\gamma^{1/2}\big(\sum_{e\in\mathcal E_h^1}\hat h_e^{-1}\|\kappa_e^{1/2}\jump{w}\|_{0,e}^2\big)^{1/2}\le\|w\|_{\EG,\kappa}$; this gives the bounds $\eta(C_{\rm tr}/\gamma)^{1/2}|||P'|||_\kappa\|w\|_{\EG,\kappa}$ and $\eta\,C_{\rm tr}^{1/2}|||P'|||_\kappa\|w\|_{\EG,\kappa}$.
For $\mathcal R_\tau$, let $r:=\nabla\cdot\kappa\nabla P'+f$, $\mathcal Dv:=\nabla\cdot\kappa\nabla v$, and $\mathcal Mv:=\kappa\nabla\omega\cdot\nabla v$.
Then $\nabla\cdot\kappa'\nabla v=(1-\omega)\mathcal Dv-\mathcal Mv$ on each element, and a direct computation gives
\[
\mathcal R_\tau(w)=\tau(\omega r,\mathcal Dw)_{\mathcal T_h}+\tau\big((1-\omega)(\omega\mathcal DP'+\mathcal MP'),\mathcal Dw\big)_{\mathcal T_h}+\tau(r-\omega\mathcal DP'-\mathcal MP',\mathcal Mw)_{\mathcal T_h}.
\]
Since $|1-\omega|\le1+\eta$, $\tau^{1/2}\|\mathcal Dw\|_{\mathcal T_h}\le\|w\|_{\EG,\kappa}$, and $\tau^{1/2}\|\mathcal Mv\|_{\mathcal T_h}\le\eta_1\|\kappa^{1/2}\nabla v\|_{\mathcal T_h}$, we obtain $|\mathcal R_\tau(w)|\le(\eta+\eta_1)(2+\eta+\eta_1)|||P'|||_\kappa\|w\|_{\EG,\kappa}$.
For $\mathcal R_N$, let $s:=\kappa\nabla P'\cdot n+j_N$.
Then $\mathcal R_N(w)=\gamma_N\langle\hat h_e(\omega s+(1-\omega)\omega\,\kappa\nabla P'\cdot n),\kappa\nabla w\cdot n\rangle_{\mathcal E_h^N}$, and $\gamma_N^{1/2}\big(\sum_{e\in\mathcal E_h^N}\hat h_e\|\kappa\nabla w\cdot n\|_{0,e}^2\big)^{1/2}\le\|w\|_{\EG,\kappa}$ gives $|\mathcal R_N(w)|\le\eta\gamma_N^{1/2}(2+\eta)|||P'|||_\kappa\|w\|_{\EG,\kappa}$.
Collecting the bounds proves~\eqref{sm:residual}.
Finally, let $\chi:=P_h(\kappa,f)-P'$.
Then $A_1^\kappa(\chi,w)=F_1^\kappa(w)-A_1^\kappa(P',w)$, and taking $w=\chi$ in~\eqref{sm:residual} together with $A_1^\kappa(\chi,\chi)=\|\chi\|_{\EG,\kappa}^2$ completes the proof.
\end{proof}

The consistency terms of $A_1$ contain traces of neural gradients, which is why $L$ depends on $n$ through $C_{\rm tr}(n,h)$, as in the continuity estimate of $A_1$~\cite[(3.30)]{kim2026locally}.
Moreover, an estimate in terms of $\eta$ alone cannot be obtained by this argument, since $\|1-\kappa'/\kappa\|_{L^\infty(\Omega)}$ can be small while $\nabla(\kappa'/\kappa)$ is not.

For inputs $u=(\kappa,f)$ with an encoded source, the source part is estimated in a discrete dual norm, which is the counterpart of $H^{-1}_D(\Omega)$.

\begin{lemma}[Discrete stability of the ENN solution]\label{sm:lem-dstab}
For admissible $(\kappa,f)$ and $(\kappa',f')$, it holds that
\[
\|P_h(\kappa,f)-P_h(\kappa',f')\|_{\EG,\kappa}\le L(\eta,\eta_1)\,|||P_h(\kappa',f)|||_\kappa+c_{\kappa,\kappa'}\,\|\ell_{f-f'}^{\kappa'}\|_{*,\kappa'},
\]
where
\[
\begin{gathered}
\ell_g^{\kappa'}(w):=(g,w-\tau\nabla\cdot\kappa'\nabla w)_{\mathcal T_h},\qquad
\|\ell\|_{*,\kappa'}:=\sup_{w\in V_h}\frac{\ell(w)}{\|w\|_{\EG,\kappa'}},\\
c_{\kappa,\kappa'}:=\sup_{0\ne v\in V_h}\frac{\|v\|_{\EG,\kappa}}{\|v\|_{\EG,\kappa'}}
\end{gathered}
\]
are the source part of $F_1^{\kappa'}$ in~\eqref{eq:A1}, the discrete dual norm, and the norm equivalence constant, respectively.
\end{lemma}

\begin{proof}
We split the difference as in the proof of \cref{lem:darcy}.
The conductivity part $P_h(\kappa,f)-P_h(\kappa',f)$ is bounded by \cref{lem:perturb}.
The source part $e'=P_h(\kappa',f)-P_h(\kappa',f')$ satisfies $A_1^{\kappa'}(e',w)=\ell_{f-f'}^{\kappa'}(w)$, so that $\|e'\|_{\EG,\kappa'}\le\|\ell_{f-f'}^{\kappa'}\|_{*,\kappa'}$.
Finally, $\|e'\|_{\EG,\kappa}\le c_{\kappa,\kappa'}\|e'\|_{\EG,\kappa'}$ by definition.
This completes the proof.
\end{proof}

Note that the first two terms of the EG norm are linear in $\kappa$, so that $c_{\kappa,\kappa'}\le e^{\varepsilon/2}$ if $\tau=\gamma_N=0$.
In general, the $\tau$-term involves $\nabla\kappa$ and the $\gamma_N$-term is quadratic in $\kappa$, so that $c_{\kappa,\kappa'}$ depends also on $\|\nabla(\log\kappa-\log\kappa')\|_{L^\infty(\Omega)}$, $\tau$, and $e^{\varepsilon}$.
With $\kappa'=\kappa_q$ and $f'=f_q$, \cref{sm:lem-dstab} is the discrete counterpart of the continuous bound~\eqref{eq:darcy-model}.

\section{Additional numerical results}\label{sm:sec-numerics}

\Cref{sm:tab-rates} lists the fitted convergence rates of \cref{sec:num-rates}.

\begin{table}[!htb]
\centering
\small
\caption{Fitted convergence rates of the network error against the exponents of \cref{thm:rate}.
The rates are least squares fits of the logarithm of the error over $N\in[64,1024]$ for the trunk ($[64,576]$ for $q=4$) and $M\in[64,2048]$ for the branch.
In the $H^1$ row, the trunk readout is the discrete $H^1$ projection, as in \cref{lem:trunk} with $t=1$.}\label{sm:tab-rates}
\begin{tabular}{lccc}
\toprule
Term, norm, input & Predicted rate & Quasi-uniform & Random (5 seeds) \\
\midrule
trunk, $L^2$, $q=2$ & $\frac12+\frac{2k_\tr+1}{2d}=2.25$ & 2.27 & 2.28 \\
trunk, $L^2$, $q=4$ & $2.25$ & 2.19 & 2.24 \\
trunk, $H^1$, $q=2$ & $\frac12+\frac{2(k_\tr-1)+1}{2d}=1.75$ & 1.78 & 1.80 \\
branch, $L^2$, $\alpha=1.5$ & $(\alpha+\frac12)/q=1.0$ & 1.10 & 1.07 \\
branch, $L^2$, $\alpha=0.5$ & $(\alpha+\frac12)/q=0.5$ & 0.59 & 0.56 \\
\bottomrule
\end{tabular}

\end{table}

\Cref{sm:tab-Clayer} lists the constants of \cref{thm:fixed-mesh-layer,cor:cons} discussed in \cref{sec:num-cons}.
Here, $C_{\rm layer}$ is computed as the square root of the largest eigenvalue of $\|S_\kappa v\|_{\EG,\kappa}^2$ relative to $\|v\|_{\EG,\kappa}^2$ on $V_n$.
Note that the numerical dimension of $V_n$ grows only from $265$ to $270$ between $n=576$ and $n=1024$, since the features become numerically dependent.

\begin{table}[!htb]
\centering
\small
\caption{Constants of the conservation layer on the fixed $8\times8$ partition for $n$ $\mathrm{ReLU}^3$ neurons with $\delta=n^{-1/2}$: the numerical dimension of $V_n$, the layer constant $C_{\rm layer}$ of \cref{lem:layer} for $\gamma_N=1$ and $\gamma_N=0$, the continuity constant $C_a$, and the trace constant $C_{\rm tr}$ of~\eqref{eq:Ctr}, maximized over $\kappa_0=1$ and $20$ random conductivities of \cref{sec:num-cons}, and the norm of the projection $\Pi$ of \cref{cor:cons} in the EG norm of $\kappa_0$.}\label{sm:tab-Clayer}
\begin{tabular}{rccccccc}
\toprule
$n$ & $\delta$ & $\dim V_n$ & $C_{\rm layer}$ & $C_{\rm layer}$, $\gamma_N=0$ & $C_a$ & $\|\Pi\|_{\rm ref}$ & $C_{\rm tr}$ \\
\midrule
36 & $1/6$ & 28 & 1.055 & 1.058 & 1.055 & 2.04 & 2.50 \\
64 & $1/8$ & 44 & 1.065 & 1.069 & 1.065 & 2.04 & 2.63 \\
144 & $1/12$ & 84 & 1.059 & 1.060 & 1.060 & 2.06 & 2.37 \\
256 & $1/16$ & 144 & 1.058 & 1.060 & 1.069 & 2.08 & 2.18 \\
576 & $1/24$ & 265 & 1.056 & 1.056 & 1.065 & 2.09 & 2.13 \\
1024 & $1/32$ & 270 & 1.055 & 1.055 & 1.061 & 2.10 & 2.02 \\
\bottomrule
\end{tabular}

\end{table}

We verified \cref{lem:perturb} for $\kappa'=\kappa_q$ with the encoders of \cref{sec:num-cons} and $50$ test conductivities, and, for the first $10$ of them, for the smaller perturbations $\kappa'=\kappa\exp(\vartheta(\log\kappa_q-\log\kappa))$ with $q=16$ and $\vartheta\in\{0.01,0.1,0.5,1\}$.
The ratio of $\|P_h(\kappa)-P_h(\kappa')\|_{\EG,\kappa}$ to the bound $L(\eta,\eta_1)|||P_h(\kappa')|||_\kappa$ is at most $0.034$, and for $\vartheta\le0.1$ the error is proportional to $\eta$, as $L(\eta,\eta_1)$ is.

\begin{figure}[!htb]
\centering
\includegraphics[width=0.66\textwidth]{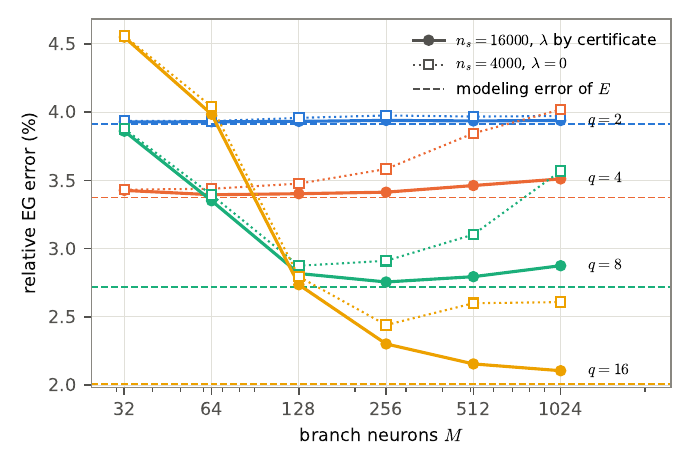}
\caption{Relative EG error of the conservative DetONet (residual loss) versus the number $M$ of branch neurons for $q\in\{2,4,8,16\}$. Solid: $n_s=16000$ with $\lambda$ selected by the certificate; dotted: $n_s=4000$ without regularization; dashed: the modeling error of $E$ for each $q$.}\label{sm:fig-cdet}
\end{figure}

\Cref{sm:fig-cdet} shows the errors of \cref{tab:cdet} versus the number $M$ of branch neurons.
The error with $n_s=16000$ stays close to the modeling error of $E$ for large $M$, whereas for $q\ge4$ the error with $n_s=4000$ and $\lambda=0$ increases for large $M$, by up to $0.7$ percentage points above its minimum over $M$, since large branches overfit the modeling error.

\bibliographystyle{siamplain}
\bibliography{references}
\end{document}